\documentclass{amsart}
\usepackage{color}

\usepackage{amsmath} 
\usepackage{amssymb}
\usepackage{mathrsfs}
 
\usepackage{mathptmx}
\usepackage{makecell}
\usepackage{tikz}

\usepackage{amsmath,amssymb,tabularx}
\usepackage{mathrsfs}
\usepackage{enumitem}

\usepackage[hidelinks]{hyperref}

\newtheorem{theorem}{Theorem}[section] 
\newtheorem{claim}[theorem]{Claim}

\newtheorem{observation}[theorem]{Observation}

\theoremstyle{definition}
\newtheorem{definition}[theorem]{Definition}

\newtheorem{convention}[theorem]{Convention}

\newtheorem{conjecture}[theorem]{Conjecture}
\newtheorem{discussion}[theorem]{Discussion}

\theoremstyle{remark}
\newtheorem{remark}[theorem]{Remark}
\newtheorem{question}[theorem]{Question}
\newtheorem{notation}[theorem]{Notation}

\newcommand{\full}{{\rm full}}
\newcommand{\almost}{{\rm almost}}
\newcommand{\medium}{{\rm medium}}
\newcommand{\semi}{{\rm semi}}
\newcommand{\pseudo}{{\rm pseudo}}

\newcommand{\ep}{{\rm ep}}

\newcommand{\AD}{{\rm AD}}

\newcommand{\tfg}{{\rm tfg}}

\newcommand{\ceq}{{\rm ceq}}

\newcommand{\ord}{{\rm ord}}

\newcommand{\PT}{{\rm PT}}

\newcommand{\CXGC}{{\rm CXGC}}

\newcommand{\ASGC}{{\rm ASGC}}

\newcommand{\BFGC}{{\rm BFGC}}
\newcommand{\AFGC}{{\rm AFGC}}
\newcommand{\CFGC}{{\rm CFGC}}

\newcommand{\SGC}{{\rm SGC}}

\newcommand{\FCG}{{\rm FGC}}

\newcommand{\WSGC}{{\rm WSGC}}
\newcommand{\CSGC}{{\rm CSGC}}

\newcommand{\PGC}{{\rm PGC}}
\newcommand{\XGC}{{\rm XGC}}
\newcommand{\JEP}{{\rm JEP}}

\newcommand{\FGC}{{\rm FGC}}

\newcommand{\WXGC}{{\rm WXGC}}

\newcommand{\AXGC}{{\rm AXGC}}

\newcommand{\AGC}{{\rm AGC}}
\newcommand{\MGC}{{\rm MGC}}

\newcommand{\pp}{{\rm pp}}
\newcommand{\qf}{{\rm qf}}

\newcommand{\TP}{{\rm TP}}

\newcommand{\otp}{{\rm otp}}
\newcommand{\feq}{{\rm feq}}

\newcommand{\suc}{{\rm suc}}

\newcommand{\NTP}{{\rm NTP}}
\newcommand{\NSOP}{{\rm NSOP}}
\newcommand{\Univ}{{\rm Univ}}
\newcommand{\club}{{\rm club}}

\newcommand{\pair}{{\rm pr}}
\newcommand{\ns}{{\rm NS}}
\newcommand{\acc}{{\rm acc}}  

\newcommand{\univ}{{\rm univ}}

\newcommand{\ZFC}{{\rm ZFC}}
\newcommand{\SOP}{{\rm SOP}}

\newcommand{\cf}{{\rm cf}}
\newcommand{\pr}{{\rm pr}}

\newcommand{\Dom}{{\rm Dom}}

\newcommand{\dmrk}{{\rm DEAR MARK:}}

\newcommand{\snn}{{\sn}}

\newcommand{\GCH}{{\rm GCH}}

\newcommand{\at}{{\rm at}}

\newcommand{\bd}{{\rm bd}}

\newcommand{\lqq}{{``}}
\newcommand{\cl}{{\rm c {\ell} }}

\def\emptycs{}
\def\evaluateROMANlist{%
        \ifx\ROMANlist\emptycs\else
        \expandafter\xxevaluate\ROMANlist\xxfertig\evaluateROMANlist\fi}
\def\xxevaluate#1,#2\xxfertig{\expandafter\newcommand\csname#1\endcsname{\mathrm{#1}}\def\ROMANlist{#2}}

 \newcommand{\ROMANlist}{Prop,pseudoProd,xyz,arey,}
 \evaluateROMANlist

\newcommand{\bfc}{{\mathbf c}}

\newcommand{\bfV}{{\mathbf V}}
\newcommand{\bfU}{{\mathbf U}}

\newcommand{\EC}{{\rm EC}}

\newcommand{\Min}{{\rm Min}}

\newcommand{\Rang}{{\rm Rang}}

\newcommand{\rest}{{\restriction}}

\newcommand{\wilog}{{\rm without loss of generality}}
\newcommand{\Wilog}{{\rm Without loss of generality}}

\newcommand{\then}{{\underline{then}}}
\newcommand{\when}{{\underline{when}}}
\newcommand{\where}{{\underline{where}}}

\newcommand{\Then}{{\underline{Then}}}

\newcommand{\Iff}{{\underline{iff}}}
\newcommand{\mn}{{\medskip\noindent}}
\newcommand{\sn}{{\smallskip\noindent}}

\newcommand{\cA}{{\mathscr A}}

\newcommand{\gb}{{\mathfrak b}}
\newcommand{\gB}{{\mathfrak B}}

\newcommand{\cE}{{\mathscr E}}
\newcommand{\varp}{{\varepsilon}}
\newcommand{\gd}{{\mathfrak d\/}} 

\newcommand{\cH}{{\mathscr H}}

\newcommand{\cF}{{\mathscr F}}
\newcommand{\cG}{{\mathscr G}}

\newcommand{\bbL}{{\mathbb L}}

\newcommand{\bbP}{{\mathbb P}}

\newcommand{\bbQ}{{\mathbb Q}}
\newcommand{\cP}{{\mathscr P}}

\newcommand{\cS}{{\mathscr S}}

\newcommand{\cX}{{\mathscr X}}

\newcount\skewfactor
\def\mathunderaccent#1#2 {\let\theaccent#1\skewfactor#2
\mathpalette\putaccentunder}
\def\putaccentunder#1#2{\oalign{$#1#2$\crcr\hidewidth
\vbox to.2ex{\hbox{$#1\skew\skewfactor\theaccent{}$}\vss}\hidewidth}}
\def\name{\mathunderaccent\tilde-3 }

\newbox\noforkbox \newdimen\forklinewidth
\setbox0\hbox{$\textstyle\bigcup$}
\setbox1\hbox to \wd0{\hfil\vrule width \forklinewidth depth \dp0
                        height \ht0 \hfil}
\setbox\noforkbox\hbox{\box1\box0\relax}
\def\unionstick{\mathop{\copy\noforkbox}\limits}
\def\nonfork#1#2_#3{#1\unionstick_{\textstyle #3}#2}
\def\nonforkin#1#2_#3^#4{#1\unionstick_{\textstyle #3}^{\textstyle
    #4}#2}
\setbox0\hbox{$\textstyle\bigcup$}
\setbox1\hbox to \wd0{\hfil{\sl /\/}\hfil}
\setbox2\hbox to \wd0{\hfil\vrule height \ht0 depth \dp0 width
                                \forklinewidth\hfil}
\newbox\doesforkbox
\setbox\doesforkbox\hbox{\box1\box0\relax}
\def\nunionstick{\mathop{\copy\doesforkbox}\limits}

\def\fork#1#2_#3{#1\nunionstick_{\textstyle #3}#2}
\def\forkin#1#2_#3^#4{#1\nunionstick_{\textstyle #3}^{\textstyle
    #4}#2}

\newcommand{\stickT}{%
\setbox255=\hbox{\raise1ex\hbox{$\hspace{0.2pt}\,\bullet\,$}}
\mathord{\rlap{\hbox to\wd255{\hss\hbox{$|$}\hss}}
\box255}
}
\newcommand{\stickS}{%
\setbox255=\hbox{\raise0.6ex\hbox{$\scriptstyle\bullet$}}
\mathord{\rlap{\hbox to\wd255{\hss\hbox{$\scriptstyle|$}\hss}}
\box255}
}

\newenvironment{PROOF}[2][\proofname.]
   {\begin{proof}[#1]\renewcommand{\qedsymbol}{\textsquare$_{\rm #2}$}}
   {\end{proof}}

\begin{document}

\title {Universality; new criterion for non-existence \\
}
\author {Saharon Shelah}  
\address{Einstein Institute of Mathematics\\
Edmond J. Safra Campus, Givat Ram\\
The Hebrew University of Jerusalem\\
Jerusalem, 9190401, Israel\\
 and \\
 Department of Mathematics\\
 Hill Center - Busch Campus \\ 
 Rutgers, The State University of New Jersey \\
 110 Frelinghuysen Road \\
 Piscataway, NJ 08854-8019 USA}
\email{shelah@math.huji.ac.il}
\urladdr{http://shelah.logic.at}
\thanks{The author thanks Alice Leonhardt for typing 
an earlier version, and an individual who wishes to remain anonymous for generously funding typing services.
We thank the ISF, The Israel Science Foundation grant
1838/19  for partially supporting this research.
This is paper number 1164 in the author's list.
First typed/compiled December 12, 2018.  References like
\cite[Th0.2=Ly5]{Sh:950} means the label of Th.0.2 is y5.
}






\date{2022-02-25}

\begin{abstract}

We find new ``reasons" for a class of models for not having a
universal model in a cardinal $\lambda$.  This work, though has
consequences in model theory, is really 
    in 
combinatorial (set theory).    
We concentrate on a
prototypical class which is a simply defined class of models, of
combinatorial character - models of
$T_{\ceq}$ (essentially another representation of 
$T_{\feq}$ which was already considered
but the proof with $T_{\ceq}$ is more transparent).
Models of $T_{\ceq}$ consist essentially of an equivalence relation on
one set and a family of choice functions for it.
This class is not simple (in the model theoretic sense) 
but seems to be very low among the non-simple
(first order complete countable) ones.  We give sufficient conditions 
for the non-existence of a universal model for it in $\lambda$. 
This work may be continued in \cite{Sh:F2071}.

\medskip 

\end{abstract}

\maketitle
\numberwithin{equation}{section}
\setcounter{section}{-1}
\newpage

\centerline{Annotated Content}
\bigskip

\noindent
\S0 \quad Introduction, (labels z,w), pg.\pageref{0}
\bigskip

\noindent
\S1 \quad For Mahlo cardinals, (label a), pg.\pageref{1}
\bigskip

\noindent
\S2 \quad On Successor Cardinals and Club Guessing, (label b), pg.\pageref{2}
\newpage

\section {Introduction} \label{0}
\bigskip

On a recent survey on the universality spectrum see \cite{Sh:1151}, an
earlier survey is \cite{Dj05}; there have been several advances
meanwhile (and this is one of the advances after \cite{Sh:1151}). See also \cite{Sh:1162}, noting the example there works also for $\mu^{+} < 2^{\mu}$ whenever $\mu$ is a strong limit singular. The problem for general first order theories is a model theoretic one,
but specific examples are combinatorial set theoretic ones (and serve
as proto-types for suitable families of theories); so
combinatorialists may ignore model theoretic notions like ``$T$ is simple,
has the tree property, is $\TP_2$", and consider only the concrete
universal theories considered; so ignore \ref{a11} (1),(2) and their
proof. Here we concentrate on the theory $T_{\ceq}$, which we considered as a 
proto-typical  {`}{`}minimal" non-simple $T$, so are expecting it (under
$\le_{\univ}$) to be low,
so is  it    
    (like $ T_{\feq}$, see below),  
        $\NSOP_1$, see
\cite{Sh:457}, \cite{Sh:692}, \cite{Sh:844}, 
\cite{MR3580894},
\cite{23-RamseyKaplan2017},
\cite{Kruckman2017}).
  True, there were non-existence
results near a strong limit singular cardinal (see on the 
$ T_{\feq}$ in \cite{Sh:457}, generalizing it 
the  oak property
\cite{Sh:710}, \cite[\S3]{Sh:820}),  
    but there were weak consistency
results on existence (see \cite{Sh:457}, \cite{Sh:614}).  We had
considered $T_{\feq}$, a prototypical example of such theories, now
$T_{\ceq}$ is essentially equivalent to it for our aims, see
\ref{a11}(3),(4)
but $ T_{\ceq} $ seem more transparent; 
we intend to deal with ``to
what family of $T$'s versions of our proof apply, in particular,
$\NTP_2$ and non-simple" elsewhere. 

We have hoped/expected that for the $\lambda > \mu = \mu^{<\mu}$ 
but $\lambda = \mu^+ < 2^\mu$ we shall have consistency results for theories like
$T_{\feq}$ and the class of triangle 
free graphs, 
\cite{Sh:1185} and hopefull  \cite{Sh:1200}.    

We first give a case with stronger
set theoretic assumptions, but more transparent proof in \S1.  In \S2
we give such proof under reasonable set theoretic assumptions, (close
to the so called club guessing) but
then have to consider finer points in combinatorial set theory on
guessing clubs.  Elsewhere we  hope to 
    have 
relevant complimentary  
consistency (see 
\cite{Sh:1185}) 
and families of theories.

A priory we think that $T_{\tfg}$, the theory of triangle free graphs,
is ``more complicated" then $T_{\feq},T_{\ceq}$, but now have doubts.

We thank Mark Po\'or and an anonymous referee for doing much to improve the paper.  

\bigskip

\begin{question}
\label{x5}
1) Does \S1 apply to more theories than in \S2?

\noindent
2) Can we characterize the dividing line?  Simple/non-simple in our context.

\noindent
3) Does it help to have:
\mn
\begin{enumerate}
\item[$(*)$]   for some $\mu,\mu < \lambda < 2^\mu$ there is
no $\cA \subseteq [\lambda]^\lambda$ which is $\mu-\AD$ of cardinality
$> \lambda$?
\end{enumerate}
\mn
This would justify the use of $\mu-\AD$ family $\cA \subseteq
[\lambda]^\lambda$ in some consistency results, see \cite{Sh:175a},
\cite{Sh:1185}, see below. 
\end{question}

\begin{discussion}
\label{x8}
Note that:
\mn
\begin{enumerate}
\item[$\boxplus$]  if $ 2 < n \le \omega,\theta \le \mu \le 
        \lambda < 
  2^\theta,\lambda \nrightarrow [\mu]^{<n}_ \theta $ 
  and we let $T_n$ be the
theory ``$\{P_k$ is an irreflexive asymmetric $k$-place relation": 
$k <  n, k \ge 2 
    \}$ and $T_n$ has a universal model $M_*$ in $\lambda$ \then \,
  there is a $\mu$-disjoint $\cA \subseteq [\lambda]^\lambda$ of
  cardinality $2^\theta$.
\end{enumerate}
\mn
[Why? Without loss of generality the universe of $M_{*}$ is $\lambda.$  Let $\bfc:[\lambda]^{< n} \rightarrow \theta  $ 
    witness $\lambda 
\nrightarrow [\mu]^{<n}_\theta $ 
    and for $u \subseteq \theta$ let $M_u =
(\lambda,\dots 
P^{M_u}_k,\ldots)_{k \in [2, n)}$   
    where $P^{M_u}_k = \{\eta \in {}^k \lambda:
\eta$ is with no repetitions and $\bfc(\Rang(\eta)) \in u\}$.
So there is an embedding $f_u$  of $ M_u $   
    into $M_*$; now $\langle   A_ u =  
    \{\pair(\alpha, f_u(\alpha ): \alpha < \lambda  \}   
:u
\subseteq \theta\rangle$ is a family as promised
when $ \pair $ is a pairing function on $ \lambda $.
    Why? If 
    $ A_{u_1}\cap A_{u_2}$ has cardinality $ \ge \mu $ 
and $ u_1 \not=  u_2 $ then 
    (letting $ B = \{ \alpha < \lambda : f_{u_1}(\alpha )= f_{u_2}(\alpha )\} $)  
\wilog \, 
$ u_1 \nsubseteq  u_2 $ and $ \mathbf{c} \rest B $  omits any 
member of $ u_1 \setminus u_2$. The rest is left to the reader.]   
%
\end{discussion}
\bigskip

\subsection {Preliminaries} \label{0B}\
\bigskip 

\begin{notation}
\label{z2}
1) $T$ is a theory with vocabulary $\tau_T = \tau(T)$ and is a
first order, if not said otherwise.

\noindent
2) 
\sn
\begin{enumerate}
\item[(a)]  $\EC_T = \{M:M$ a model of $T\}$,
\sn
\item[(b)]  $\EC_T(\lambda) = \{M \in \EC_T:M$ of cardinality $\lambda\}$,
\sn
\item[(c)]  $\EC_T(\lambda!) = \{M \in \EC_T:M$ has universe
  $\lambda\}$, 
\item[(d)]  for a set $ A $ of ordinals and ordinal
$ \alpha $ let $ \suc_A ( \alpha ) $ be 
$ \min \{ \beta  \in A : \beta > \alpha \} $.  
\end{enumerate}   

\noindent 
3) Let $ \pair$ be  be an (easily computable)  pairing function on  
ordinals such that 
for $ \alpha, \beta $   we have
$\pair (\alpha,  \beta ) < \max \{\omega, \alpha + | \alpha |, \beta + |\beta  | \} $. 
\end{notation}

\begin{convention}
\label{z6}
1) 
\mn
\begin{enumerate}
\item[(A)]  If $T$ is a f.o. 
theory not complete
(like $T^0_{\ceq},T^0_{\feq},$ usually universal), \then \, embedding are the usual ones,
(on 
 $\EC_T$) and $\subseteq_T$ 
 (on $ \EC _T$) 
 means $\subseteq$
  and we assume $\EC_T$ has
amalgamation and $\JEP$.
\sn
\item[(B)]  If $T$ is complete, \then \, embeddings are elementary 
(on $\EC_T$)
  and $\subseteq_T$ means $\prec$ on $\EC_T$.
\sn
\item[(C)]  We say $f$ is a $T$-embedding of $M$ into $N$ or $f:M
  \longrightarrow_T N$ \when \, $M,N$ are models of $T,f$ embed $M$
  into $N$ and $f(M) \subseteq_T N$.
\end{enumerate}
\mn

1A) In any case we always assume $T$ has $\JEP$ (for $\subseteq_{T}$ of course). 

\mn

2) If $\Delta \subseteq \bbL(\tau_T)$ (such that $T$ has $\JEP$ under $\Delta-$embedding) then
$\univ_{T,\Delta}(\lambda)$ is the minimal $\chi$ such that there
is a sequence $\bar M$ which is a $(\lambda,T,\Delta)$-universal
sequence which means: 
\mn
\begin{enumerate}
\item[(a)]  $\bar M = \langle M_\alpha:\alpha < \chi\rangle$ is a sequence 
of models of $T$, 
\sn
\item[(b)]  each $M_\alpha$ is of cardinality $\lambda$,  
\sn
\item[(c)]  for every model $M$ of $T$ of cardinality $\lambda$ there
  is a $\Delta$-embedding of $M$ into some $M_\alpha$, see below.
\end{enumerate}
\mn
3) For given $T,\Delta$ as above and models $M,N$ of $T$, we say $f$ is
a $\Delta$-embedding of $M$ into $N$ \when \,:
\mn
\begin{enumerate}
\item[(a)]  $f$ is a function from $M$ into $N$, 
\sn
\item[(b)]  if $\varphi(x_0,\dotsc,x_{n-1}) \in \Delta$ and
  $a_0,\dotsc,a_{n-1} \in M$ and $M \models
  \varphi[a_0,\dotsc,a_{n-1}]$ \then \, $N \models
  \varphi[f(a_0),\dotsc,f(a_{n-1})]$, 
\sn
\item[(c)]  so $f$ is one-to-one when $(x \ne y) \in \Delta$.
\end{enumerate}
\mn
4) For $T,\Delta$ as above in part (2) we may omit $\Delta$ \when \,:
\mn
\begin{enumerate}
\item[(a)]  $T$ is complete, $\Delta = \bbL(\tau_T)$, all first order
formulas, 
\sn
\item[(b)]  $T$ not complete, $\Delta$ the set of quantifier  free 
    formulas
  in $\bbL(\tau_T)$.
\end{enumerate}
\mn
5) We may write $\at,\ep$ for $\Delta_{\at(T)} = \{\varphi \in
  \bbL(\tau_T):\varphi$ is atomic$\},\Delta_{\ep(T)} = \{\varphi \in
  \bbL(\tau_T):
  \varphi$ existential positive$\}$ respectively.  We may
  write $\tau$ instead of $T$.  We may write $\varphi$ instead $\Delta
  = \{\varphi\}$ and $\pm \varphi$ instead $\Delta = \{\varphi,\neg \varphi\}$.
\end{convention}

\begin{notation}
\label{z9}
1) Let $\alpha,\beta,\gamma,\delta,\varp,\zeta,\xi,i,j$ denote
ordinals.

\noindent
2) Let $\kappa,\lambda,\mu,\chi,\partial,\theta,\Upsilon$ denote
cardinals, infinite if not said otherwise.

\noindent
3) Let $k,\ell,m,n$ denote natural numbers.

\noindent
4) Let $\varphi,\psi,\vartheta$ denote formulas, 
f.o. if not said otherwise.

\noindent
5) For $\lambda > \kappa$ regular cardinals let $S_{\kappa}^{\lambda} = \{ \delta < \lambda: \cf(\delta) = \cf(\kappa)\}$ and $S_{\leq \kappa}^{\lambda} = \{ \delta < \lambda: \cf(\delta) \leq \kappa \}.$
\end{notation}

\begin{definition}
\label{z12}
1) $J^{\bd}_\theta = \{A \subseteq \theta:\sup(A) < \theta\}$,
$\bd$ stands for bounding, for $\theta$ a regular cardinal 
or just a limit ordinal.

\noindent
1A) For $\theta$ regular uncountable let:
\mn
\begin{itemize}
\item  $D^{\club}_\theta = \{A \subseteq \theta$: there is a club (=
closed unbounded subset) $E$ of $\theta$ such that $E \subseteq A\}$.
\item $ \ns_ \theta $ is the non-stationary ideal on $ \theta $.  
\end{itemize}
\mn
2) For a regular $\theta$ let:
\mn
\begin{enumerate}
\item[(a)]  $\gd_\theta = \Min\{|\cF|:\cF \subseteq {}^\theta \theta$
  is $<_{J^{\bd}_\theta}$-cofinal in $\mu_\mu\}$
\sn
\item[(b)]  $\gb_\theta = \Min\{|\cF|:\cF \subseteq {}^\theta \theta$
has no $<_{J^{\bd}_\theta}$-upper bound$\}$.
\end{enumerate}
\mn
3) Let $\gd^{\club}_\theta$ be defined similarly using
$<_{\ns  
    _{\theta}}$ when $\theta$ is regular uncountable.
    
\noindent 
4)
For a model $ M $ and a set $ u \subseteq M $ let $ M \rest u $  
is defined naturally, allowing a function symbol  to be interpreted 
as a partial function (and so an individual constant   to be not defined) but $M \rest u \subseteq M$ means $u = \cl_{M}(u),$ see below.

\noindent
5) For a model $M$ and $A\subseteq M$ let $\cl_{M}(A) = \cl(A, M)$ be the minimal subset $B$ of $M$ including $A$ and closed under the functions of $M;$ so $M \rest \cl_{M}(A) \subseteq M$ and if $M$ has Skolem functions then $M \rest \cl_{M}(A) \prec M.$ 

\end{definition} 

Recall   
\label{88r-0.5}
\begin{definition} \label{z16}   
1) For a regular uncountable 
cardinal $\lambda$ let $\check I[\lambda] = \{S \subseteq \lambda$: some
pair $(E,\bar a)$ witnesses $S \in \check I(\lambda)$, see below$\}$. \newline  

\noindent 
2) We say that $(E,u)$ is a witness for $S \in \check I[\lambda]$
\underline{if}:
\begin{enumerate}
\item[{$(a)$}]  $E$ is a club of the regular cardinal $\lambda$, 
\sn
\item[{$(b)$}]  $u = \langle u_\alpha:\alpha < \lambda
\rangle,u_\alpha \subseteq \alpha$ and $\beta \in u_\alpha \Rightarrow
u_\beta = \beta \cap u_\alpha$, 
\sn
\item[{$(c)$}]  for every $\delta \in E \cap S,u_\delta$ is an
unbounded subset of $\delta$ of order-type $< \delta$ (and $\delta$ is
a limit ordinal, necessarily $ \delta $ is not a regular cardinal).  
\end{enumerate}

\noindent   
3) For $ \kappa =  \mathrm{cf}(\kappa) < \lambda = \cf( \lambda ) $  let 
$ \check{I}_  
    {\le \kappa} [ \lambda ]$  be the ideal 
    $ \{ S \subseteq \lambda : S \subseteq S^\lambda _{\le \kappa}, 
        S \in \check I[\lambda] \} $  
\end{definition}

\label{88r-0.6}  
By  (\cite{Sh:108}, \cite{Sh:88a} and better) \cite{Sh:420} and \cite{Sh:E12} 
we have: 
\begin{claim}   \label{z19}
Let $\lambda$ be regular uncountable. \newline 
1) If $S \in \check I[\lambda]$ {\then} 
we can find a witness $(E,\bar a)$ for $S \in \check I[\lambda]$ such that 
(clauses (a), (b), (c) from \ref{z16}(2) and):  
\begin{enumerate}
\item[{$(d ) $}]  $\delta \in S \cap E \Rightarrow \,{\text{\rm otp\/}}
(a_\delta) = \,{\text{\rm cf\/}}(\delta)$,  
\sn
\item[{$(e)$}]  if $\alpha \notin S$ then ${\text{\rm otp\/}}
(a_\alpha) < \,{\text{\rm cf\/}}(\delta)$ for some $\delta \in S \cap E$.
\end{enumerate}
2) $S \in \check I[\lambda]$ {\Iff} there is a pair $(E,\bar{\mathcal P})$
such that:
\begin{enumerate}
\item[{$(a)$}]  $E$ is a club of the regular uncountable $\lambda$, 
\sn
\item[{$(b)$}]  $\bar{\mathcal P} = \langle {\mathcal P}_\alpha:\alpha <
\lambda\rangle$, where ${\mathcal P}_\alpha \subseteq \{u:u \subseteq
\alpha\}$ has cardinality $< \lambda$,  
\sn
\item[{$(c)$}]  if $\alpha < \beta < \lambda$ and $\alpha \in u \in 
{\mathcal P}_\beta$ then $u \cap \alpha \in {\mathcal P}_\alpha$,  
\sn
\item[{$(d)$}]  if $\delta \in E \cap S$ then some $u \in 
{\mathcal P}_\delta$ is an unbounded subset of $\delta$ 
    of order type 
        $ < \delta $  
(and $\delta$ is a limit
ordinal).
\end{enumerate}
\end{claim}


\noindent
3) We say a stationary subset $S$ has club guessing when some $\langle C_{\delta}: \delta \in S \rangle$ witnesses it, which means: $C_{\delta}$ is a club of $S$ and for every club $E$ of $\lambda$ for some $\delta \in S$ we have $C_{\delta} \subseteq \lambda.$   

\section {On $T_{\ceq}$ for Mahlo cardinals} \label{1}
 \bigskip

As here we consider $T_{\ceq}$ the simplest, non-simple
theory, we may consider how much does it behave like the class of
graphs (equivalently random graph)?  We prove that
not by a non-existence result, but with quite specific set
theoretic assumptions.

$ T_{\ceq}$ is very close to (and equivalent for our purposes to) 
the older 
$T_{\feq}$ which   
        is 
a prime example for a theory with the tree order
property,   equivalently non-simple   
(even $\TP_2$ but 
    having 
neither the strict order property nor even
just the $\SOP_1$).
For it we get  here parallel and better results than \cite{Sh:457} where 
it is proved that there are
limitations on the universality spectrum for $T_{\feq}$ and in
\cite{Sh:710},  
    which
generalize the results 
for any $T$ with the 
so called oak property, 
see somewhat more in \cite[\S3]{Sh:820}.  The results 
in those papers 
are 
meaningful when SCH fails, that is, consider a cardinal $\lambda$ such that:
for some strong limit singular $\mu,\mu^+ < \lambda < 2^\mu$   
if $\lambda$ is regular   
then ``usually" $T_{\feq} $,   
has no universal in $\lambda$. 
 
But what about $\lambda \in (\mu,2^\mu)$ when
for transparency we assume $\mu = \mu^{< \mu}$?.  
    Here (in \S1)  
    we get further such non-existence results for 
        (weakly inaccessible)  
    Mahlo cardinals.
In \S2, we do better
but the Mahlo case may cover more classes,
comes first and the proofs are more transparent.  The proof here  
    (in \S1)  
can be axiomatized as in \S2 using:
\mn
\begin{enumerate}
\item[$\boxplus$]  $\PGC(\lambda,S)$  where $ S $ is   
    a stationary set of regular cardinals $ < \lambda $  
    means that some  
    $ \mathbf{U} $ witness $\PGC(\lambda,S)$  where  
    $ \mathbf{U} = \{ \langle \omega ( 1 + \varepsilon) : \varepsilon < \theta 
        \rangle: \theta \in S \} $
(so $  \mathbf{U}  
        =S$).
        
        \noindent
        Recall that ``$\mathbf{U}$ witness $\PGC(\lambda, \theta)$''  means $\PGC(\lambda, \theta) = \min \{ \vert \mathbf{U} \vert : \mathbf{U} \subseteq \mathbf{U}_{\lambda, \theta} \ \text{and} \ \mathbf{U}  \\ \text{does} \ $P-$ \text{guess clubs} \}.$ See \ref{b2}(5),  Definition \ref{b2}(3c) and \ref{b2}(4). 
\end{enumerate}
\bigskip

\noindent
First, recall (the reader can concentrate on the universal versions,
$T^0_{\feq},T^0_{\ceq}$,   on $T_{\feq}$ see \cite[2.1=Lb3,3.1=Lc3]{Sh:457}): 
\begin{definition}  
\label{a2}
$T_{\feq} = T^1_{\feq}$ is the model completion of the
following   (universal first order) 
theory, $T^0_{\feq}$ which is defined by:
\mn
\begin{enumerate}
\item[(A)]  $\tau = \tau(T^0_{\feq})$ consists of:
\sn
\begin{enumerate}
\item[(a)]  predicates $P,Q$ (unary),  
\sn
\item[(b)]  $E$ (three place predicate written as 
    $x E_z y $ instead $ E(x,y,z)$),  
\end{enumerate}
\sn
\item[(B)]  a $\tau$-model $M$ is a model of $T^0_{\feq}$ \Iff \,:  
\sn
\begin{enumerate}
\item[(a)]  the universe of $M$ is the disjoint union of $P^M$ and $Q^M$,
\sn
\item[(b)]   
  $ x E_z y\rightarrow P(z) \wedge Q^M(x) \wedge Q^M(y)$,  
\sn
\item[(c)]  for any fixed  
    $ z\in P^M,E^M_z$ is an equivalence relation on
  $Q^M$.
\end{enumerate}
\end{enumerate}
\end{definition}

\begin{observation}
\label{a5}
0) $ T_{\feq}$ is well defined and    
    $ \univ( \lambda , T_{\feq} )=   
    \univ (\lambda  ,  T^0_{\feq}) $   

\noindent 
1) So if $M \models T_{\feq}$ \then \,:
\mn
\begin{enumerate}
\item[$(*)$]
\begin{enumerate}
\item[(a)]  in (B)(c)  of Def. \ref{a2},  
for each $x \in P^{M}, E_{x}^{M}$ 
is with infinitely many equivalence classes,  
\sn
\item[(b)]  if $n < \omega,x_1,\dotsc,x_n \in P^{M}$ with no repetition
  and $y_1,\dotsc,y_n \in Q^{M}$ then for some $y \in
  Q^{M},\bigwedge\limits^{n}_{\ell =1} y E^{M} _{x_\ell} y_\ell$,
\sn
\item[(c)]  if $n < \omega$ and $y_1,\dotsc,y_n \in Q^{M}$ and $e$ is an
  equivalence relation on $\{1,\dotsc,n\}$ \then \, for some $x \in P^{M}$
  we have $y_\ell E^{M} _x y_k \Leftrightarrow \ell e k$,
  \snn 
\item[(d)]    $ P^M, Q^M $ are infinite.  
\end{enumerate}
\end{enumerate}
\mn
2) Hence $T_{\feq}$ has elimination of quantifiers 
    and $ \univ_{T_{\feq}}(\lambda )= \univ_{T^0_{\feq}}(\lambda )$.  
\end{observation}  

\noindent
We present a close relative, the main one we consider here (and, as
proved below, equivalent to $T_{\feq}$ for our purpose).
\begin{definition}  
\label{a8}    
\mn
    $T_{\ceq} = T^1_{\feq}$ is the model completion of the
following   (universal first order) 
theory, $T^0_{\ceq}$ which is defined by:  
\begin{enumerate}
\item[(A)]  $\tau = \tau(T^0_{\ceq}) = \tau(T_{\ceq})$ consists of: 
$P,Q$ unary predicates, $E$ a binary predicate and $F$ a binary function symbol,
\sn
\item[(B)]  a $\tau$-model $M$ is a model of the universal theory.  
$T^0_{\ceq}$ \Iff \,:  
\sn
\begin{enumerate}
\item[(a)]  $P^M,Q^M$ is a partition of $M$,
\sn
\item[(b)]  $E^M$ is an equivalence relation on $Q^M$,
\sn
\item[(c)]  $F^M$ is a function from $Q^M \times P^M$ into $Q^M$ such
that for every $c \in P^M, a \mapsto 
    F^M(a,c)$ is choosing a representative for the 
$a/E^M$-equivalence class, that is, we have:
\sn
\begin{enumerate}
\item[$(\alpha)$]  $a \in Q^M \Rightarrow F^M(a,c) \in a/E^M$,
\sn
\item[$(\beta)$]  if $a,b \in Q^M$ are $E^M$-equivalent  
then
  $F^M(a,c) = F^M(b,c)$. 
\item[$(\gamma)$]  if $ c \notin P^M \vee a  \notin Q^M $  
    then $ F^M(a,c)$ is not defined (or, if you prefer, is equal to $c$).  
\end{enumerate}
\end{enumerate}
\end{enumerate}
\end{definition} 
\mn
Concerning $\lambda$ in the neighborhood of a strong limit singular we
shall not give details as we can just quote.

\begin{claim}  
\label{a11}
0) Concerning $ T^0_{\ceq}$
\begin{enumerate} 
\item[(a)] For a model $ M $ of $ T^0_{\ceq}$ and  
 $ A \subseteq M $ with $ n $ elements, the closure of $ A $
 inside $ M $ has at most $ n + n^2$ elements, (even at most 
 $ n + (n/2)^2$ elements),
 \item[(b)]  $ T^0_{\ceq}$ has amalgamation and JEP,  
 \item[(c)] $ T^0_{\ceq}$ has a model completion, that is $ T_{\ceq}$
 is well defined, 
 \item[(d)]  $ \univ( \lambda, T_{\ceq})  
    = \univ( \lambda , T_{\feq})   
    = \univ_{T^0_{\ceq}}(\lambda )$.  
\end{enumerate} 

\noindent 
1) $T_{\ceq}$ is not simple, is $\NSOP_2$ and even $\NSOP_1$ and has the
oak property, in fact, by $\qf$ (quantifier free) and even atomic formulas.

\noindent
2) We have $(A) \Rightarrow (B)$ \where \,:
\mn
\begin{enumerate}
\item[(A)]
\begin{enumerate}
\item[(a)]   $\theta < \mu < \lambda < \chi$,
\sn
\item[(b)]  $\cf(\lambda) = \lambda,\theta = \cf(\theta) =
  \cf(\mu),\mu^+ < \lambda$,
\sn
\item[(c)]  $\chi := \pp_{\Gamma(\theta)}(\mu) > \lambda + |i^*|$,
\sn
\item[(d)]   there is a set $\{(a_i, b_i):i<i^*\}$ with $a_i \in [\lambda]^{<\mu},
b_i \in [\lambda]^\theta$ and $|\{b_i:i< i^*\}| \le \lambda$ such that:
for every $f:\lambda \rightarrow \lambda$ for some $i,f(b_i) \subseteq a_i$,
\end{enumerate}
\sn
\item[(B)]
\begin{enumerate}
\item[(a)]   $T_{\ceq}$ equivalently $T^0_{\ceq}$ 
has no universal model in $\lambda$,
\sn
\item[(b)]   Moreover, $\univ(\lambda,T_{\ceq}) \ge \chi = 
\pp_{\Gamma(\theta)}(\mu)$.
\end{enumerate}
\end{enumerate}
\mn
3) $T_{\feq}$ can be interpreted in $T_{\ceq}$ hence
$\univ_{T_{\feq}}(\lambda) \le \univ_{T_{\ceq}}(\lambda)$.

\noindent
4) Also the inverse of part (3)  
    holds.
\end{claim}

\begin{PROOF}{\ref{a11}}
0) Easy as clause (d) follows by parts (3), (4).

\noindent
1) By (1) (and part (3), (4)) quoting \cite{Sh:710} where the oak property was
introduced.

\noindent 
2) Follows from parts (3), (4) and 
    \cite[Claim 2.2]{Sh:457}. 
    
\noindent
3) For a model $M$ of $T^0_{\ceq}$ we define a model $N=N[M]$ of
$T^ 0  
_{\feq}$ as follows:  
\mn
\begin{enumerate}
\item[$(*)_{N,M}$]
\begin{enumerate}
\item[(a)]  $Q^N = P^M,P^N = Q^M/E^M$,  
\sn
\item[(b)]  $E^N = \{(a,b,C):C \in P^N= Q^M/E^M$ and $a,b \in Q^N$ and 
$(\forall c \in C)[F^M(c,a) = F^M(c,b)]$ equivalently, 
$(\exists c \in C)[F^M(c,a) = F^M(c,b)]\}$.  
\end{enumerate}
\end{enumerate}
\mn
Now check 
that $ N \models T^0  
    _{\feq}$ and $ M \models T_{ceq} \Leftrightarrow N \models T_{feq} $. 

\noindent
4) For a model $N$ of $T^0_{\feq}$ we define a model $M= M[N]$ of
$T^0_{\ceq}$ as follows:
\mn
\begin{enumerate}
\item[$(*)_{M,N}$]
\begin{enumerate}
\item[(a)]   $ P^M= Q^N $  
    and $Q^M = \{(c, A): c \in P^N $ and   
    $ A \in 
    Q^N/E^N_c \}   $ 
\sn
\item[(b)] $ E^M = \{ ((c_1, A_1), (c_2, A_2)):
c_1 = c_2 \in P^M $ and  
$ A_1, A_2 \in Q^N/E^N_{c_2}\} $
\sn 
\item[(c)]  $F^M:Q^M \times P^M \rightarrow Q^M$ is defined by: If \,
  $d \in Q^M,b \in P^M$ hence for some $c \in P^N,A \in   
        Q^N /  
        E^N_c $ we 
  have $d = (c,A)$ \then \, we let $F^M(d,b) = (c, 
  b/E^M_c). $ 
\end{enumerate}
\end{enumerate}
\mn
Now check
that $ N \models T^0  
        _{\feq}$ and $ M \models T_{\ceq} \Leftrightarrow N \models T_{\feq} $. 
\end{PROOF}  





\noindent
We now point out a new reason involved ``large $\gd_\theta$'s" for  
not having a universal model in $\lambda$,
even for many non-simple $T$'s.  In this section we deal with a case where the
proof is simpler using $T_{\ceq}$ and $\lambda$ a Mahlo cardinal.
\begin{claim}
\label{a14}
    1)
Assume $\lambda$ is a (weakly inaccessible) Mahlo cardinal and $S =
\{\theta < \lambda:\theta$ regular (weakly inaccessible) and 
$\gd  
    _\theta > \lambda\}$ is stationary in $\lambda$  
and $S$ has club guessing.  

\Then
\mn
\begin{enumerate}
\item[(a)]  $\univ(\lambda,T_{\ceq})$ is $> \lambda$,  
\sn
\item[(b)]  even,  
    $\ge \sup\{\chi^+  $:    
the set $\{\theta \in
 S:\gd_\theta > \chi\}$ is stationary
    and has club guessing$\}$. 
        
\end{enumerate}

\noindent 
2) We have $ \chi < 
\univ( \lambda, T_{\ceq} ) $    
\when:
\begin{enumerate} 
    \item[(a)] $ \lambda $ is a Mahlo  weakly inaccessible cardinal,
    \item[(b)] $ \lambda  \le 
            \chi $, 
    \item[(c)] $ S \subseteq \{\theta  < \lambda :  
        \theta $ is weakly inaccessible  
cardinal$\} $ 
is stationary.  
    \item[(d)] $ \bar{ {\mathscr P} } = 
\langle  {\mathscr P} _ \theta  : \theta  \in S\rangle $, 
    \item[(e)]  if $ \theta \in S $  then  $  {\mathscr P} _ \theta $ is a set of $ \le \lambda $  
clubs of $ \theta  $,
    \item[(f)] $ \bar{ {\mathscr P} } $ guess clubs of $ \lambda $, that is,
for every club  $ E $ of $ \lambda $ for some 
$ C \in {\mathscr P}_ \theta , \theta  \in S$ 
we have $ C \subseteq E$,
    \item[(g)] $ \mathfrak{d} _ \theta  > 
    \chi $   
        for every $ \theta \in S $. 
\end{enumerate} 
\end{claim}


\begin{PROOF}   {\ref{a14}}
1)  
Clearly
\mn
\begin{enumerate}
\item[$(*)_0$]  
     it suffices to:
  \begin{enumerate} 
        \item[(a)]     
        fix $\chi \ge \lambda$ such that  
  $S_\chi = \{\theta \in S:  \mathfrak{d}   
  _\theta >   
  \chi \} $  
  is stationary and has club guessing,  

        \item[(b)] prove $\univ_T(\lambda) >\chi$.
\end{enumerate}  
\end{enumerate} 
\mn
Let $T=T_{\ceq},$ without less of generality assume $S = S_{\chi}$ and let:
\mn
\begin{enumerate}
\item[$(*)_1$]  $\langle C_\delta:\delta \in S\rangle$ witness
  ``$S$ has club-guessing"; 
\sn
\item[$(*)_2$]  if (A) below holds, \then \, we define some objects in
  (B) where:
\sn
\begin{enumerate}
\item[(A)]
\begin{enumerate}
\item[(a)]   $M \in \EC_T(\lambda!)$, 
\sn
\item[(b)]  $|P^M| = \lambda  $,\
\snn 
\item[(c)] 
    $ \theta$ is regular and    
     $ \theta \in S $, 
\sn
\item[(d)]   $E$ a club of $\theta$. 
\end{enumerate}
\sn
\item[(B)]  we define:
\sn
\begin{enumerate}
    \item[(a)]  for $a \in  P^M $ 
hence  $a < \lambda$ let $g_a = g_{M,E,a}$ be
  the following  function from $\theta$ to $\theta$: 
\sn
\begin{itemize}
    \item  for $\alpha < \theta,g_a(\alpha)$ is   
the minimal $\beta \in E$
such that: $\beta \in E \backslash (\alpha +1)$ and $(\beta_1
\in Q  
    ^M \cap \beta) \wedge (F^M(\beta_1,a) < \theta) \Rightarrow
F^M(\beta_1,a) < \beta$, 
\end{itemize}
\sn
\item[(b)]  $\cG^0_{M,E} = \{g_{M,E,a}:a \in P  
            ^M\}$, note that $ E $ determine $ \theta $,  
\sn
\item[(c)]   for $ \theta \in S $ let  
    $\cG^*  
        _{M,\theta} = \{ g_{M,C_\theta,a}:a \in P  
                ^M \}$.
\end{enumerate}
\end{enumerate}
\end{enumerate}
\mn
Now easily
\mn
\begin{enumerate}
\item[$(*)_3$]  for $a,M,\theta,E$ as above, $g_{M,E,a}$ is a well defined non-decreasing function from $\theta$ into $\theta$, in fact, into $E
  \subseteq \theta,$
\sn
\item[$(*)_4$]  if $M,N \in \EC_T(\lambda!)$ and $f$ embeds $M$ into $N$
\then \, for some club $E^*$ of $\lambda$: if $\theta \in S,\theta =
\sup(E^* \cap \theta),E \subseteq E^* \cap \theta$ is a club of $\theta$
and $a \in P^M$ then $g_{M,E,a} \le g_{N,E,f(a)} $ (so, the only way $E^{*}$ influences is the demand ``$E \subseteq E^{*}$''). 

\end{enumerate}
\mn
Recall that, 
$\theta \in S \Rightarrow \gd_\theta > \chi \ge
\lambda$  
 and we shall prove that $\univ(\lambda,T) > \chi$; this
suffices.  So assume $\langle M_\alpha:\alpha < \chi\rangle$ is 
a sequence of members of
$\EC_T(\lambda !)$. 

So for each $\theta \in S$ the set $\cG_\theta = 
\cup\{\cG^*  
    _{M_\alpha,\theta}:\alpha < \chi $ satisfies  
    $ | Q^{M_\alpha} \cap \theta| = \theta\}$ has cardinality $\le \chi$   
recalling $\lambda \le \chi$.

For $ \theta \in S $, 
as $|\cG  
        _\theta| < \gd 
            _\theta$, 
    necessarily there is   an increasing 
$g_\theta \in {}^\theta \theta$ such that $g \in \cG_\theta
\Rightarrow   
 g_ \theta  \nleq  g \mod J^{\bd }_ \theta $ and without loss of generality, $ g_{\theta} \in {{}^{\theta} ( C_{\theta}}).$  
Now we define a model $N \in \EC_{T^0_{\ceq}}  
            (\lambda!)$ with $\tau_N =
\tau(T_{\ceq})   = \tau(T^0_{\ceq}) $ as follows:
\mn
\begin{enumerate}
\item[(A)]   universe is $\lambda$, 
\sn
\item[(B)]
\begin{enumerate}
\item[(a)]  $ Q^N $  
        is the set of odd ordinals $< \lambda$, 

\sn
\item[(b)]  $ E^N $   
    is an equivalence relation on $Q  
        ^N$ such that for
  every $\alpha < \beta < \lambda$ satisfying $\beta$ is divisible by
  $|\alpha|,\alpha \in Q  
        ^N$ we have $|\alpha/E^N \cap \beta| = |\beta|$, 

\sn
\item[(c)]   if $\alpha = 4 \beta +1 < \lambda$ then $\alpha/E^N$ is
  disjoint to $\alpha$, 
        
\sn
\item[(d)]   if $\theta \in S$ and $ \alpha < \theta $
    then $  \theta >  
            F^{N}(4 \alpha + 1 ,  
        \theta) >
  g_\theta(4 \alpha  + 1)$.  
\end{enumerate}
\end{enumerate}
\mn
This is easy to do. 

To show that $ \bar{ M } $ does not witness   
$ \univ( \lambda , T) \le \chi $  is suffice to show that $ N $ 
cannot be embedded in $ M_ \alpha $ for any $ \alpha < \chi $.  
Toward contradiction assume that $ \alpha < \chi $ and $ f $ 
is an embedding of $ N $ into $ M_ \alpha $.    
Let $ E = \{\delta < \lambda : \delta $ a limit ordinal such that 
    $ ((M_ \alpha \rest \delta , N \rest \delta , f \rest \delta, <  \rest \delta  ) \prec 
        (M_ \alpha , N, f , < \rest \lambda ))\}  $.  
        Clearly $ E $ is a club of $\lambda $ hence 
        for some $ \theta \in S_ \chi $ 
        we have $ C_ \theta \subseteq E $. Let 
        $ h \in {}^{ \theta }\theta  $ be 
        $g_{M_ \alpha, C_ \theta,   
            f( \theta )}$, so 
        is well defined and belongs 
        to $ {\mathscr G} ^0 _{M_ \alpha , C_ \theta }$   
         hence to  
          ${\mathscr G} ^*_{M_ \alpha , \theta }$ 
           hence to $ {\mathscr G} _ \theta $  
        hence $ g _ \theta \not < h \mod J^{\bd}_ \theta.$ Now,





\begin{enumerate}
    \item[$\bullet_{1}$] choose $\alpha < \theta$ such that $h(\alpha) < g_{\theta}(\alpha),$ 
    
    \sn
    
    \item[$\bullet_{2}$] let $\gamma = 4 \alpha + 1,$
    
    \sn
    
    \item[$\bullet_{3}$] $F^{N}(\gamma, \theta) \in (g_{\theta}(\alpha), \theta)$ by the choice of $N$ i.e. (B)(d) above,
    
    \sn
    
    \item[$\bullet_{4}$] $g_{\theta}(\alpha) \in C_{\theta}$ by the choice of $g_{\theta},$
    
    \sn
    
    \item[$\bullet_{5}$] $C_{\theta} \subseteq E$ by the choice of $\theta,$
    
    \sn
    
    \item[$\bullet_{6}$] $g_{\theta}(\alpha) \in E$ by $\bullet_{4} $ and $ \bullet_{5},$
    
    \sn
    
    \item[$\bullet_{7}$] every member of $E$ is closed under $f$ and $f^{-1},$ 
\end{enumerate}
    
[Why? By the choice of $E$ it is closed under $f,$ as $f$ is  one-to-one similarly for $f^{-1}$.]
    
\begin{enumerate}
    \item[$\bullet_{8}$] $f(F^{N}(\gamma, \theta)) \in [g_{\theta}(\alpha), \theta),$
\end{enumerate}    
    
[Why? By $\bullet_{3},\bullet_{6} $ and $ \bullet_{7}.$]
    

\begin{enumerate}
     \item[$\bullet_{9}$] $f(F^{N}(\gamma, \theta)) = F^{M_{\alpha}}(f(\gamma), f(\theta)),$
\end{enumerate}
    
[Why? As $f$ embed $N$ into $M_{\alpha}.$]

\begin{enumerate}

    \item[$\bullet_{10}$] $F^{M_{\alpha}}(f(\gamma), f(\theta)) \in [g_{\theta}(\alpha), \theta),$ 

\end{enumerate}    

[Why? by $\bullet_{8} $ and $ \bullet_{9}.$]

\begin{enumerate}
    \item[$\bullet_{11}$] $h(\alpha)$ is a  member of $C_{\theta},$ hence a limit ordinal and in $E,$
\end{enumerate}    
 
[Why? By the choice of $h.$]

\begin{enumerate}    
    \item[$\bullet_{12}$] $\alpha < h(\alpha),$ $\gamma < h(\alpha) $ and $ f(\gamma) < h(\alpha),$
\end{enumerate}

[Why? First, $\alpha < h(\alpha)$ by the choice of $h.$ Second, $\gamma < h(\alpha)$: as $h(\alpha)$ is limit $> \alpha$  by $\bullet_{11}$ and $\gamma = 4\alpha +1$ by $\bullet_{1}.$ Third $f(\gamma) < h(\alpha)$: as $\gamma < h(\alpha)$ and $h(\alpha) \in E$ by $\bullet_{11}$ and so $h(\alpha)$ is closed under $f$ by $\bullet_{7}.$] 
    

\begin{enumerate}    
    \item[$\bullet_{13}$] $F^{M_{\alpha}}(f(\gamma), f(\theta)) < h(\alpha),$
\end{enumerate}

[Why? By the choice of $h$ as $g_{M_{\alpha}, C_{\theta}, f(\theta)}$ and $\bullet_{12}$].

Now by the inequalities $\bullet_{1}, \bullet_{8}, \bullet_{9}$ and $ \bullet_{13}$ we get
$h(\alpha) < g_{\theta}(\alpha) \leq f(F^{N}(\gamma, \theta)) = F^{M_{\alpha}}(f(\gamma), f(\theta)) < h(\alpha),$ contradiction.

\noindent 
2) Similarly.  
\end{PROOF}

\begin{remark}   \label{a19}  

Under the assumption of \ref{a14}(2), we can similarly prove that:
for every sequence $ \langle (E_ \xi , {\mathscr G} _{ \xi, \theta }  )   
:  
\xi < \chi, \theta \in S \rangle  $  
satisfying clause (A) below, there is a sequence 
$ \langle g_{\theta, \alpha }: \theta \in S,\alpha < \lambda \rangle $
with $ g_ {\theta, \alpha }\in {}^{ \theta } \theta $
satisfying clause (B)   
    below, where:
\begin{enumerate} 
    \item[(A)]  $ E_ \xi $  is a club of $ \lambda$ for $ \xi < \chi $
        and $ {\mathscr G} _{\xi, \theta }\subseteq {}^{ \theta } \theta $ 
        has cardinality $ \le \lambda $ for $ \xi < \lambda, \theta \in S$.  
    \item[(B)] for every $ \xi < \chi $ and club $ E $ of $ \lambda $
there are $ \theta \in \acc(E_ \xi ) \cap E \cap S $ 
and $ \alpha < \lambda $ such that $\theta = \sup(E_{\xi} \cap E \cap \theta)$ and
$ g \in {\mathscr G} _{\xi , \theta }  \Rightarrow g_{\theta , \alpha }
    \nleq g \mod J^{\bd} _{E_ \xi \cap E \cap \theta }$.
\end{enumerate} 
\end{remark} 
\newpage

\section {On successor Cardinals and club guessing} \label{2}
\bigskip

We first introduce the relevant notions (in \ref{b2}); (we could add clause
\ref{b2}(2)(b) into the definition of $\bfU_{\lambda,\theta}$ in  
\ref{b2}(1), 
but so far it does not matter\footnote{Note 
    that it is relevant to 
``fully $ D $-guess clubs" implies ``almost  
guess clubs", see \ref{b5}}).  
We then investigate it and use it for
sufficient conditions for \lqq no universal".  

\begin{definition}
\label{b2}
Assume $\lambda > \theta $ 
    are regular and $D \subseteq
\cP(\theta)$ 
    is upward closed   non-empty 
satisfying $D \subseteq [\theta]^\theta$,  
omitting $D$ means $D = \{\theta\}$; 
and $ \mathfrak{B} $ is a model with universe $ \lambda $ 
and countable vocabulary 
    but $ \mathfrak{B} $ is locally finite when $ \theta = 
    {\aleph_0} $.  
Saying \lqq for $ D $-most 
$ \varepsilon < \theta $"  
will mean
\lqq for some $ X \in D $  for every $ \varepsilon \in X $". 
The main case\footnote{ 
    We may omit clause (b) from the definition \ref{b2}(3) 
    of \lqq fully $ D$-guess clubs",  
    the only problem this cause is for it implying the other
    versions, (see \ref{b5}).  
}
is $ \theta > {\aleph_0} $, this  
        is necessary for 
the \lqq full" cases (see parts (2)), but not for the others;  
we may forget to assume $ \theta > {\aleph_0}$.  

\noindent
1) Let $\bfU_{\lambda,\theta} = \{\bar u:\bar u = \langle u_i:i <
\theta\rangle$ is $\subseteq$-increasing continuous, and $i < \theta
\Rightarrow i \subseteq u_i \in [\lambda]^{<\theta}$ (hence 
$ \theta \subseteq \cup\{u_i:i < \theta\} \in [\lambda]^\theta)$
    and 
$\bigwedge\limits_{i < \theta }  
    u_i \cap \theta \in
\theta \}.$  

\noindent
1A) We shall say that $ \bfU  \subseteq \mathbf{U} 
_{\lambda, \theta } $ obeys $ \mathfrak{B} $ 
when every $ \bar{ u } \in \bfU $ does, which 
means that  for every $ \varepsilon < \theta $ 
we have $ \mathfrak{B} \upharpoonright u_ \varepsilon 
\subseteq \mathfrak{B} $, (if $ \mathfrak{B} $ 
has Skolem functions this is equivalent to
 $  \mathfrak{B} \upharpoonright u _ \varepsilon 
 \prec  \mathfrak{B} $  which implies $ \theta > {\aleph_0} $).  

\noindent
2) We say $\bfU \subseteq \bfU_{\lambda,\theta}$ fully $D$-guesses clubs \when \, 
    $ \theta >  {\aleph_0} $ and  
for every model $M$ with universe $\lambda$ and countable vocabulary
there is $\bar u \in \bfU$ which fully $D$-guesses $M$ meaning\footnote{
    We may omit  in \ref{b2}(2)  
    the clauses $ (a)( \alpha ), (b)$
    but then we have problems with $\lqq  \FGC  \Rightarrow \AGC$ 
    and the gain is doubtful.  
} 
\begin{enumerate} 
    \item[(a)]   
\begin{enumerate} 
        \item[$(\alpha )$] if\footnote{This implies a case of club guessing.} $ \varepsilon < \theta $  then  $  \cl( u_ \varepsilon , M)  \subseteq \sup(u _ \varepsilon ) $,  
        moreover (actually follows using an expansion of $ M $)
        $ M \rest \sup( u_ \varepsilon ) \prec M,$

        \item[$(\beta )$]  
$(\exists \cX
\in D)(\forall \varp)[\varp \in \cX \Rightarrow \cl_{M}(u_{\varepsilon}) = u_{\varepsilon} 
 \subseteq M]$, i.e. 
    for $ D $-most $ \varepsilon < \theta $  the set 
$u_\varp$ is closed under the 
functions of $M$, (in an equivalent definition 
$M_\varp \rest u_\varp
\prec M$ as we can expand $M$ by Skolem functions).
        \end{enumerate} 
     \item[(b)]  the sequence $ \ord( \bar{ u }) = 
        \langle \sup( u_ \varepsilon) :
      \varepsilon < \theta \rangle $        
is strictly increasing.  
  \end{enumerate} 


\noindent
3) We say $\bfU \subseteq \bfU_{\lambda,\theta}$ almost $D$-guesses
clubs 
\when \,:
\mn
\begin{enumerate}
\item[(a)] 
for every model $M$ with universe $\lambda$ and countable
  vocabulary and $A \in [\lambda]^\lambda$ for some $\bar u \in \bfU$ 
we have: 
\sn
\begin{enumerate} 
\item[$(\alpha)$] if $\varp < \theta$ then $\cl(u_{\varp}, M) \subseteq \sup(u_{\varp})$;  as in $(a)( \alpha) $ of part (2) without the moreover,
\sn
\item[$(\beta)$]  
        for $ D $-most $ \varepsilon < \theta $  we have   
  $ A \cap u_{\varp +1} \nsubseteq \sup(u_\varp)$, 
\sn
\item[$(\gamma)$]  $c \ell(\bigcup\limits_{\varp < \theta 
    } u_\varp,M) =
\bigcup\limits_{\varp  < \theta   
        } u_\varp$, that is, $M \rest
(\bigcup\limits_{\varp <  \theta 
    } u_\varp) \subseteq M$, 
\end{enumerate}
\sn
\item[(b)]  if $\bar u \in \bfU$ then $\ord(\bar u) = 
\langle \sup(u_\varp):\varp < \theta\rangle$ is strictly increasing.
\end{enumerate}
\mn
    3A) We say $\bfU$ medium $D$-guesses clubs \when \, as in 
      part 
(3) omitting  clause 
$(a)(\gamma)$.

\noindent
3B) We say $\bfU \subseteq \bfU_{\lambda,\theta}$ semi-$D $-guesses clubs  
\when \,:
\mn
\begin{enumerate}
\item[(a)$'$]  as  (a)  
        in part (3) but replacing $(\beta)$ by:
\sn
\begin{enumerate}
\item[$(\beta)'$]  for $ D $-most   
$\varp < \theta$   
for some $\zeta \in  
  [\varp,\theta)$ and $\alpha \in A$ we\footnote{
    The \lqq $ \alpha \notin  u _\zeta$'' follows, and ``$D$ -most'' can be replaced by ``all''.
  }
  have $\alpha \in
(u_{\zeta +1} \backslash u_\zeta) \cap 
(\sup(u_{\varepsilon  +1}) \backslash
\sup(u_\varp))$,  
\end{enumerate}
\sn
\item[(b)]  as in part (3). 
\end{enumerate}
\mn
3C) 
    We say $\bfU \subseteq \bfU_{\lambda,\theta}$ pseudo-$ D $-guess clubs
\when \,:
\mn
\begin{enumerate}
\item[(a)$''$]  if 
    $ M $ is as above and 
$A \in [\lambda]^\lambda$ then for some $\bar u \in
\bfU$ we have:  
    \begin{enumerate} 
        \item[$( \alpha ) $]  as is part (3) 
        clause $(a)(\alpha )$,  
        \item[$(\beta )$]   for $ D $-most  
        $\varp < \theta$ 
            for some $\zeta \in  
  [\varp,\theta)$ and $\alpha \in A$ we have $\alpha \in
(u_{\zeta +1} \backslash u_\zeta) \cap 
(\sup(u_{\varepsilon  +1}) \backslash
\sup(u_\varp))$,
    \end{enumerate} 
\sn
\item[(b)]  as above.
\end{enumerate}
\mn
3D) We say $\bfU$ is $(\lambda,\theta)$-reasonable 
(or just reasonable   
    when $ ( \lambda , \theta ) $ are clear
from the context)
\when \, $\bfU \subseteq 
\bfU_{\lambda,\theta}$ satisfies clause (3)(b).



\noindent
4) We say $\bfU$ does $X-D$-guess clubs 
 \when \,:
\mn
\begin{itemize}
\item  $\bfU $ does fully
$D$-guess clubs and $X=F$, 
\sn
\item  $\bfU$ does almost $D$-guess clubs and $X=A$, 
\sn
\item  $\bfU$ does semi-guess clubs and $X=S$, 
\sn
\item  $\bfU$ does medium $D$-guess clubs and $X = M$, 
\sn
\item  $\bfU$ does pseudo guess clubs and $X=P$.
\end{itemize}
\mn
5) Let $\XGC_D(\lambda,\theta) = \min\{|\bfU|:\bfU \subseteq
\bfU_{\lambda,\theta}$ and $\bfU$ does $X-D$-guess clubs$\}$. 

\noindent
5A) Similarly  $ \XGC _D (\lambda , \theta , \mathfrak{B} )$
when we restrict ourselves to $ \mathbf{U} $ 
obeying $ \mathfrak{B} $.

\noindent
6) We say $\bfU \subseteq \bfU_{\lambda,\theta}$ is bounded \when \,
there is an $F$ witnessing it which means: $F$ is a function from 
$\{\bar u \rest ( \zeta + 1 ) : \bar u \in \bfU,  
\zeta < \theta  
\}$ 
into $\lambda$ such that $F(\bar u_1 \rest ( \zeta_1 + 1 )  
) = F(\bar
u_2 \rest (\zeta_2 + 1 ) 
) \Rightarrow \bar u_1 \rest (\zeta_1 + 1 ) 
= \bar u_2 \rest
(\zeta_2 + 1) $  
and $F(\bar u \rest (\zeta +1)) < \sup(u_{\zeta +1})$.

\noindent
7) We say ``strongly bounded" when in addition $F(\bar u \rest
(\zeta +1)) \in u_{\zeta +1}$ for every $\zeta < \theta$.


\noindent
8) We say $\bfU \subseteq \bfU_{\lambda,\theta}$ is weakly bounded,
\when \, there is a function $F$ witnessing it which means:
\mn
\begin{enumerate}
\item[(a)]  $\Dom(F) = \{\ord(\bar u \rest (\zeta + 1)) 
    :\bar u \in \bfU$ and
  $\zeta < \theta  
    \}$ recalling 
  $\ord(\bar u)=
  \langle \sup(u_\varp):\varp < \theta\rangle$, 
\sn
\item[(b)]  $\Rang(F) \subseteq \lambda$ and $F(\ord(\bar u) \rest
  (\zeta +1)  ) 
    < \sup(u_{\zeta +1})$ for $\bar u \in \bfU$ and 
$\zeta < \theta$, 
\sn
\item[(c)]  if $\zeta_1,\zeta_2 < \theta$ are successor of successor 
ordinals and 
$\bar u_1,\bar u_2 \in \bfU$ and $F(\ord(\bar u_1) \rest \zeta_1) 
= F(\ord(\bar u_2)\rest 
\zeta_2)$ then $\ord(\bar u_1) \rest \zeta_1 
= \ord(\bar u_2) \rest \zeta_2$.
\end{enumerate}

\noindent
9) Let
\mn
\begin{enumerate}
\item[(a)]  if $\bar u \in \bfU_{\lambda,\theta}$ and $f:\theta
\rightarrow \theta$ is $ \le $-increasing continuous 
with limit $ \theta $ 
then $\bar u^{[f]} = 
\bar u[f]:=\langle u_{f(\varp)}:\varp < \theta\rangle]$, 
\sn
\item[(b)]  if $\bfU \subseteq \bfU_{\lambda,\theta}$ and $f:\theta
  \rightarrow \theta$ is $ \le $-increasing continuous 
  with limit $ \theta $ 
  then $\bfU^{[f]} := \bfU[f] 
= \{\bar u 
        {[f]}:\bar u \in \bfU\}$, 
\sn
\item[(c)]  if $\bfU \subseteq \bfU_{\lambda,\theta}$ 
and $ \cF $ is a set of $ \le $-increasing 
continuous function from $ \theta $ into
$ \theta $ with limit $ \theta $ 
  then $\bfU[\cF] = \{\bar u[f]:\bar u \in
  \cF,f \in \cF\}$, 
\sn
\item[(d)]  if $w \in 
    [\theta]^ \theta$ then $f_w = f[w]$ is the 
$g:\theta \rightarrow \theta$ such that $\langle g(\varp): \varp < \theta$ but the closure of $w$ in order. 
  so is $ \le $-increasing continuous with limit
  $ \theta $.
\end{enumerate}
\mn
10)   In (a),(b)
of part (9) 
above we may write $\bar u[w],\bfU[w]$ for $w \in
[\theta]^\theta$ meaning $\bar u[f],\bfU[f]$ where $f=f_w$, 
writing $\bfU[W],W \subseteq
[\theta]^\theta$ mean $\cup\{\bfU[w]:w \in W\}$.

\noindent 
11)
Now for $X \in \{F,A,S,M,P\}$ we let 
(naturally   
and we can add $ \mathfrak{B} $ as in part (5A)
):
\mn
\begin{enumerate}
\item[(a)]  $\AXGC_D(\lambda,\theta) = \Min\{|\bfU|:\bfU \subseteq  
  \bfU_{\lambda,\theta}$ does $X-D$-guess clubs and is strongly bounded$\}$, 
\sn
\item[(b)]  $\CXGC_D(\lambda,\theta)$ is defined as in (a) but $\bfU$ is
  just bounded, 
\sn
\item[(c)]  $\WXGC_D(\lambda,\theta)$ is defined as in clause (a) 
but $\bfU$ is weakly bounded.
\end{enumerate}
\end{definition}

\noindent
Some of the obvious implications are:
\begin{observation}
\label{b5}
1) If $\bfU$   
    fully $D$-guesses  clubs, 
\then \, $\bfU$ almost $D$-guesses clubs,    


\noindent
2) If $\bfU$ almost $D$-guesses clubs \then \, $\bfU$ 
semi-guess-club and medium $D$-guesses clubs.

\noindent
3) If $\bfU$ semi-$D$-guesses-clubs or medium $D$-guesses clubs \then \, $\bfU$ does pseudo $ D $-guesses clubs.

\noindent
4) If $D_1 \subseteq D_2 \subseteq [\theta]^\theta$ 
then ``$ \mathbf{U} $   
does $X- 
    D_1$-guess clubs"
implies \lqq$\bfU$ does  
$X-    D_2$-guess clubs" for 
$X \in  \{ F,A,M,   S,P  \}$, 
    we may write $\{\full,\almost,\medium,  \semi, \pseudo 
            \}$. 


\noindent 
5) Assume $ \mathbf{U} \subseteq \mathbf{U} _{\lambda , \theta }$   
and $ \mathfrak{B} $ is as in \ref{b5}. 
Then there is $ \mathbf{U} '  $  
such that: 

\begin{enumerate} 
\item[(a)] $   \mathbf{U}' \subseteq \mathbf{U} _{\lambda,\theta }$ 
\item[(b)] $ | \mathbf{U}' |  \le | \mathbf{U} | $ 
\item[(c)]  if $ \mathbf{U} $ does  $X-  D $-guess clubs, 
for $ X \in \{ F,A,M, S, P \} $   
as in part (4)   
then so does $ \mathbf{U} ' $,
\item[(d)] $ \mathbf{U} $ obeys  $ \mathfrak{B} $,  (see \ref{b2}(1)). 
\end{enumerate} 

\noindent 
6) In \ref{b2}(11) the number is  $ \ge \lambda $. 

\noindent 
7) We may replace \lqq countable vocabulary" by 
\lqq vocabulary of cardinality $ < \theta$''.
\end{observation}

\begin{PROOF} {\ref{b5}}
E.g.    

\noindent 
5)  
Let $ \mathbf{U} ' = \{ \bar{ u }' \in \mathbf{U} : $  
for some $ \bar{ u } \in \mathbf{U} $ for every $ \varepsilon < \theta $ 
we have $ \cl_ \mathfrak{B} ( u_\varepsilon )  = 
u'_ \varepsilon \subseteq \sup( u_ \varepsilon )  \},  $  
    it suffice to prove that 
    $ \mathbf{U} ' $ is as required. The main point is to 
    verify the appropriate version  of clause (a)  
        in Def \ref{b2}.
    So let $ M $ be a model with universe $ \lambda $ and countable
    vocabulary, we have to find a suitable 
        meber of 
        $ \mathbf{U} '$.
    By renaming, \wilog \, the vocabulary of $M $ 
    is disjoint to 
    the one of $ \mathfrak{B} $ 
    and let $ M' $ be a common expansion 
    of $ M$ and $ \mathfrak{B} $ 
            with $ \tau (M' )=  
        \tau (M) \cup \tau ( \mathfrak{B} )$. 
        Let $ E = \{ \delta : M 
        \rest \delta \prec  M \} $.    
        So $ (M', E , < \rest \lambda )$
    is as required in clause (a) 
    for $ \mathbf{U} $ hence there is a suitable  $ \bar{ u } \in \mathbf{U} $. 
    We can check that in all cases  
    $ \bar{ u }'   
    = \langle  \cl _\mathfrak{B} (u_ \varepsilon ): \varepsilon < \theta   \rangle  
    \in \mathbf{U} 
    $ is as required here, so we are done. 
    
    \noindent 
    7) 
    Recall \ref{b2}(1), the statement ``$u_{\epsilon} \cap \theta \in \theta$''. 
\end{PROOF} 

\begin{definition}
\label{b8}
1) For the model theory: for a model $M \in \EC_T(\lambda!),\Delta
\subseteq \bbL(\tau_T)$ and $u \subseteq \lambda,A \subseteq M$ 
let $M^{[A]} \rest_\Delta u$ be the model
$M \rest u$ expanded by all  the restriction  to  $u$ 
of all relations definable by a
$\Delta$-formula with 
parameters from $A$.

\noindent
1A) If $\Delta = \bbL_{\qf}(\tau_M)$ \then \, we may omit $\Delta$; 
writing $\bar a$ instead $A$ means $\Rang(\bar a)$.

\noindent
2) For $M \in \EC_T(\lambda!),\bar u \in \bfU_{\lambda,\theta}$ and
$\bar a \in {}^{\omega >}M$ let $g_{\bar a,\bar u,M}$ be the function
from $\theta$ to $\theta$ such that for $\zeta < \theta,g_{\bar
  a,\bar u,M}(\zeta)$ is the minimal $\varp \in (\zeta,\theta)$
such that $M^{[\bar a]} \rest u_\varp \prec M^{[\bar a]} \rest 
\bigcup\limits_{\xi  < \theta   
        } u_\xi$.
\end{definition}

\begin{claim}
\label{b11} 
We assume  that $ \mathfrak{B} $ is a model 
with universe $ \lambda $ and countable vocabulary, (for the case of full club guessing, we add locally finite when $ \theta={\aleph_0}.) $  

\noindent 
1) We have 
\mn
\begin{enumerate}
\item[(A)] 
$\CSGC(\lambda,\theta) = \lambda$, 
moreover $ \lambda = \CSGC ( \lambda,\theta  , \mathfrak{B})  $ 
provided that:  
    \begin{itemize}  
        \item  
            $\lambda = \cf(\lambda) =  
                \theta^{++}$ and $\theta = \cf(\theta) $ 
    \end{itemize} 
\sn
\item[(B)]  
    $\ASGC(\lambda,\theta) = \lambda$
  provided\footnote{ 
  We can weaken the demand: if we weaken the demand
    in Definition \ref{b2}(5) to ``for stationary many $\varp <
    \theta$" and $\theta \ge \aleph_2$.
    } 
    that
\sn
\begin{itemize}
  \item $ \lambda = \theta ^{++}, 
    \theta = \cf( \theta ) $, 
  \item   
    there is a stationary set 
$S \subseteq S^{\theta^+}_\theta$ from
    $\check I_\theta[\theta^+]$,  
\end{itemize}
\sn
\item[(C)] $\AFGC(\lambda,\theta)   = \lambda  $   even with a reasonable 
  witness. provided that:
\begin{itemize}  
\item
$\lambda = \lambda^\theta  $   
        and $\theta = \cf(\theta)  > {\aleph_0}  $,  
\end{itemize}  
 
\snn
\item[(D)] 
    $ \MGC_D(\lambda, \theta )= \lambda $ \when:  
\begin{enumerate} 
\item[(*)]  $ \theta = \cf( \theta )  < \lambda $ 
    and there is  $ {\mathscr S } $  such that:  
        \begin{enumerate}
    \item[(a)]        $ {\mathscr S } \subseteq \{w : w \subseteq \lambda ,
            \otp( w) = \theta \} $, 
    \item[(b)] $ {\mathscr S } $ has cardinality $ \lambda $,  
    \item[(c)] if $ A \in [\lambda ]^\lambda $ then  
        for some $ w \in {\mathscr S } $ 
            the set $ w \cap A $ has cardinality $ \theta $,  
    \item[(d)]          $ D = [ \theta ] ^ \theta $.
        \end{enumerate} 
\end{enumerate} 
  \item[(E)]  $ \AGC(\lambda, \theta )= \lambda $ \when \, 
\item[(*)] we have:   
\begin{enumerate}  
    \item[(a),(b),(c),(d)]  as in (D) 
            above, 
    \item[(e)]  the cofinality of $ ([\lambda ]^ \theta , \subseteq ) $ is equal to 
        $ \lambda $.

\end{enumerate} 
\end{enumerate}
\mn
2)  
For regular $\lambda > \theta = \cf(\theta) $  we have:  
\begin{enumerate} 
    \item[(A)]   if $ \SGC_D  (\lambda, \theta )= \lambda $ and  
        $ \mathfrak{b} _ \theta \le \lambda $  \then \, 
        $ \AGC_D ( \lambda, \theta )= \lambda $ when $ D= [\theta ]^\theta $, 
        
    \item[(B)] if $ \SGC( \lambda, \theta ) )= \lambda $ and 
        $\mathfrak{d} _ \theta   \le \lambda  $ \then \, 
        $\AGC ( \lambda, \theta ) = \lambda $
         recalling that the default
        value of $ D $  is $ \{ \theta \}$.     
\end{enumerate}

\mn
3) For $\lambda > \theta = \cf(\theta)$  
 such that\footnote{see footnote to part (2)
}    
    $ \lambda >   
            \theta ^+ $
we have
$\SGC(\lambda,\theta) = \lambda$ provided that (e.g. $\lambda =
\theta^{+n}$ for some $n > 0 $):  
\mn
\begin{enumerate}
\item[$\boxplus^3_{\lambda,\theta}$]  $\cf([\lambda]^\theta,\subseteq)
  = \lambda$.
\end{enumerate}
\mn
4) If $\bfU_1 \subseteq \bfU_{\lambda,\theta}$ medium  
  guesses  
    clubs,
\then \, there is $\bfU \subseteq \bfU_{\lambda,\theta}$ which  medium  
guesses  
    clubs of cardinality $\le |\bfU_1|$ and  for $ \bar{ u } \in \mathbf{U} $
    we have:  
\mn
\begin{enumerate}
\item[(a)]  if $ u = \cup \{ u_i: i < \theta 
\} $ 
then $u \subseteq \delta = \sup(u)$ for
  some $\delta < \lambda$ of cofinality $\theta$; 
    (this actually follows by \ref{b2}(3)(b)),  
\sn
\item[(b)]  if $\gb_\theta \le \lambda$ then 
$ u = \cup \{ u_i: i < \theta 
\} $ 
and $ \bar{ u }  \in \mathbf{U} _{\lambda, \theta }  $
then  
$\otp(u  \setminus \theta ) = \theta$, 



\end{enumerate}
\mn
5)  
If $\lambda \ge \theta^+$ and $\theta = \cf(\theta) > \aleph_0$ and
$S \subseteq \{\delta < \lambda:\cf(\delta)=\theta\}$ is stationary
and some $\bar C = \langle C_\delta:\delta \in S\rangle$ guesses  
    clubs, \then \,
$\PGC(\lambda,\theta) = \lambda$. 

\noindent 
6) 
If  $ \cf([\lambda ]^\theta . \subseteq )=   \lambda, \theta > {\aleph_0}  $ 
and $ \mathfrak{d} _ \theta \le \lambda  $ 
then 
$ \FGC ( \lambda, \theta ) 
    = \lambda $, moreover  
$ \BFGC ( \lambda, \theta ) = \lambda $, (in fact,
looking at \cite{Sh:420} we get strongly bounding).

\end{claim}

\begin{discussion}
\label{b14}  



\noindent 
1)  In \ref{b11} 
we have \ZFC \, results,
we  may 
get stronger results (on the full 
and almost versions) in some forcing extensions
see \ref{b47}. and \cite{Sh:F2071}.

\noindent 
2) We can look at the cases of Definition \ref{b2} for singular $\lambda$,
replacing $(u_\zeta \backslash \sup (u_\varp)) $ 
by $u_  \zeta 
\backslash u_\varp$,  
        but we have not arrive to it.  

\noindent
3) When we have clause (a)$(\gamma)$ of the Definition \ref{b2}(3)
there is  less 
need of clause $(a)(\alpha)$.
E.g. in \ref{b11}(1)(C)  
we do not need ``$\lambda$ regular".

\noindent 
4)
In clauses (D), (E) of \ref{b2}(1) we may add bounded/weakly bounded under
    natural assumption. 
\end{discussion}

\begin{PROOF}{\ref{b11}}  
\Wilog \, $ \mathfrak{B}  $ has a pairing function   
$  {\rm pr}^\mathfrak{B}$ and its inverses as well as $\alpha + 1, \alpha + \beta $ and $\alpha \beta.$  

\noindent 
1) \underline{Clause (A)}: First, 
 choose 
$S,S^+ , 
    \bar C$ such that (partial square guessing clubs):
\mn
\begin{enumerate}
\item[$(*)_1$]
\begin{enumerate}
\item[(a)]  $S \subseteq \{\delta < \lambda:\cf(\delta) = \theta$ and
  $\delta > \theta^+\}$ is stationary, 
\sn
\item[(b)]  $S \subseteq S^+ \subseteq \{\delta < \lambda:\cf(\delta)
\le \theta$  and $ \delta > \theta ^+ \} $, moreover
if $ \delta \in S $ then $ \delta = \sup (S^+ \cap \delta ) $ ,
\sn
\item[(c)]  $\bar C = \langle C_\alpha:\alpha \in S^+\rangle$,
\sn
\item[(d)]  $C_\alpha$ is a closed subset of $\alpha$ of order type
  $\le \theta$, and 
  $\otp(C_\alpha)$ is a limit ordinal \Iff \, 
$\alpha = \sup(C_\alpha)$,  
\sn
\item[(e)]  for $\alpha \in S^+$ we have $\alpha \in S \Leftrightarrow
  \otp(C_\alpha)=\theta$,
\sn
\item[(f)]  if $\alpha \in C_\beta$ then
        $ \alpha \in S^+ $  and 
    $C_\alpha = C_\beta \cap
  \alpha$,
\sn
\item[(g)]   $\bar C \upharpoonright S $ 
 guess clubs, i.e.: if $E$ is a club of $\lambda$
  then for stationarily many $\delta \in S$ we have $C_\delta
  \subseteq E$,
\sn
\item[(h)]  if $\alpha \in S^+$ then $ \alpha > \theta ^+ $  and 
        $\alpha$ is closed under
  $ \mathfrak{B} $, that is $ \mathfrak{B} 
  \rest \alpha \subseteq \mathfrak{B} $, 
\item[(i)]  
    if $ \alpha \in S^+$ then $ \theta ^2 $  divides $ \delta $.
\sn
\end{enumerate}
\end{enumerate}
\mn
Why do they exist (provably in $\ZFC$)? 
 see \cite[1.3=L1.3(b)]{Sh:E12}, but we elaborate 
    (for the case $ \theta > {\aleph_0} $);   
by \cite[4.4(1), pg.47]{Sh:351} (with $\theta^{+}, \lambda$ here standing for $\lambda, \lambda^{+}$ there):

\begin{enumerate}
    \item[$(*)_{1.1}$] there are $W, \overline{S}, \overline{C}_{i} \ (i < \theta^{+})$ such that: 
        \begin{enumerate}
            \item[(A)]  $W = \{  \delta < (\theta^{+})^{+}: \cf(\delta) < \theta^{+} \}$  hence is in $\check{I}[(\theta^{+})^{+}].$
            
            \item[(B)] $W$ is the union of $\lambda$ sets which have the square property, i.e., there are sequences $\overline{S} = \langle S_{i}: i < \lambda \rangle$ and $\overline{C}_{i} = \langle C_{\delta}^{i}: \delta \in S_{i} \rangle$ for $i < \lambda$ such that: 
                \begin{enumerate}
                    \item[(a)] $ W \subseteq \bigcup_{i < \lambda} S_{i},$
                    
                    \item[(b)] For $\delta \in S_{i},$ $C_{\delta}^{i}$ is a subset of $\delta \cap S_{i}$ of cardinality $< \lambda$ closed in $\delta,$ and if $\delta$ is a limit ordinal then $C_{\delta}^{i}$ is unbounded in $\delta,$ 
                    
                    \item[(c)] For all $\delta_{1}, \delta_{2}$ if $\delta_{2} \in S_{i}$ and $\delta_{1} \in C_{\delta_{2}}^{i}$ then  $\delta_{1} \in S_{i}$ and $C_{\delta_{1}}^{i} = C_{\delta_{2}}^{i}\cap \delta_{1}.$ (Notice that $\delta_{1}$ may also be a successor ordinal.)
                \end{enumerate}
        \end{enumerate}
\end{enumerate}

Easily (and as in \cite[a. III]{Sh:e} making $\theta^{+}$ tries):

\begin{enumerate}
    \item[$(*)_{1.2}$] there are $\zeta < \theta^{+}, \ i < \lambda$ and a club $E_{*}$ of $\lambda$ such that: for every club $E \subseteq E_{*}$ of $\lambda$ for some $\delta$ we have $\delta \in S_{i}, \ \cf(\delta) = \theta, \ \delta = \sup(C_{\delta}^{i} \cap E_{*}), \ C_{\zeta}^{i} \cap E = C_{\delta} \cap E_{*}$ and $\otp(C_{\delta}^{i} \cap E_{*}) = \zeta.$ 
    
    \item[$(*)_{1.3}$] without loss of generality $\alpha \in E_{*} \Rightarrow (\cl_{\mathfrak{B}}(\alpha) = \alpha) \wedge (\theta^{+})^{2} \mid  \alpha     \wedge \alpha \geq \theta^{+},$

    \item[$(*)_{1.4}$] let: 
        \begin{enumerate} 
            \item[(a)] $e \subseteq \zeta$ is unbounded in $\zeta$ and $\otp(e) = \theta$
            
            \item[(b)] $S = \{ \delta \in S_{i}: \cf(\delta) = \theta, \ \otp(C_{\delta}^{i} \cap E_{*}) = \zeta \},$
            
            \item[(c)] $S^{+} = \{ \alpha: \alpha\in S \ \text{or for some} \  \delta \in S \ \text{we have} \ \alpha \in C_{\delta}^{i} \ \text{and} \ \otp(\alpha \cap C_{\delta}^{i}) \in e \}.$
            
            \item[(d)] $C_{\delta} = \{ \alpha \in C_{\delta}^{i}: \otp(\alpha \cap C_{\delta}^{i}) \in e \}$ for $\delta \in S.$
        \end{enumerate}
\end{enumerate}

Now $S, \delta^{+}, \langle C_{\delta}: \delta \in S^{+} \rangle$ satisfies all the demands, proving $(*)_{1}.$

\begin{enumerate}
\item[$(*)_2$]  For $\delta \in S$ let 
$\langle \gamma^\bullet_{\delta,\varp}:\varp <
  \theta\rangle$ list $C_\delta$ in increasing order.
\end{enumerate}
\mn
Second, fix $\bar f,\bar g$ such that:
\mn
\begin{enumerate}
\item[$(*)_3$]
\begin{enumerate}
\item[(a)]  $\bar f = \langle f_\alpha:
\alpha \in [\theta^+,\lambda)\rangle$,
\sn
\item[(b)]  $f_\alpha$ is a one-to-one function from $\theta^+$ onto
  $\alpha$,
\sn
\item[(c)]  $\bar g = \langle g_\xi : \xi 
\in  [\theta,\theta^+)\rangle$,
\sn
\item[(d)]  $g_\xi 
$ is a one-to-one function from $\theta$ onto
  $\xi   $.
\end{enumerate}
\end{enumerate}
\mn
Third, 
\mn
\begin{enumerate}
\item[$(*)_4$]
\begin{enumerate}
\item[(a)]  for $\delta \in S$ let $e_\delta = \{
\xi  < \theta^+$: if
  $\alpha \in C_\delta$ then $\Rang(f_\alpha \rest 
  \xi ) = \alpha
  \cap \Rang(f_\delta \rest 
  \xi )$ 
  and this set
  includes $ C_ \delta  \cap \alpha   $ and 
  has cardinality $\theta\}$
\sn
\item[(b)]  $e_\delta$ is a club of $\theta^+$.
\end{enumerate}
\end{enumerate}
\mn
[Why clause (b) holds?  As $\otp(C_\delta) 
= \theta$ and $ \alpha \in C_  \delta 
    \cup \{ \delta \} 
\Rightarrow | \alpha | = \theta ^+ $,
 this should be clear.]
\mn
\begin{enumerate}
\item[$(*)_5$]   for $\delta \in S$ and $\xi \in e_\delta$ let:
\sn
\begin{enumerate}  
\item[(a)]   $u_{\delta,\xi} = \Rang(f_\delta 
\rest \xi)  $, it belongs to 
 $  [\delta]^\theta$ and it includes $ C_ \delta $,  
\sn
\item[(b)]  we choose $\bar u_{\delta,\xi} = 
\langle u_{\delta,\xi,\varp}:\varp
< \theta\rangle$ by $u_{\delta,\xi,\varp} = 
c {\ell}  _ \mathfrak{B} ( 
\{f_{\gamma^\bullet_{\delta,\upsilon}}  
(g_\xi(\zeta)):\upsilon < \omega (1+\varp)$
  and $\zeta < \omega (1+ \varepsilon ) \}
  \cup \{ \gamma ^ \bullet _{\delta , \upsilon}:
  \upsilon < \omega (1+\varepsilon )\}) $,  
\sn
\item[(c)]  for $w \in [\theta]^\theta$ let 
$\bar u^{[w]}_{\delta,\xi}$ be 
$ \langle u ^{[w]}_{\delta, \varepsilon 
}; 
  \varepsilon < \theta \rangle $
   where 
   $ u^{[w]}_{\delta, \varepsilon }  =   
 u_{\delta,\iota} $  
  where:    
    $\iota \in w$  is the minimal $ \iota  $ that   
    satisfies 
 $\otp(w
  \cap \iota)= \varepsilon $,  
  this fits \ref{b2}(9)(d).
  
\end{enumerate}
\end{enumerate}
\mn
Note that (recalling $(*)_2$):

\mn
\begin{enumerate}
\item[$(*)_6$]  For $ \delta \in S , \xi \in e_ \delta $
we have: 
\begin{enumerate}
\item[(a)]  $\bar u_{\delta,\xi}$ is a $\subseteq$-increasing continuous
  sequence of subsets of $u_{\delta,\varp}$,  
\sn
\item[(b)]  each $u_{\delta,\xi,\varp}$ 
include $ C_{\gamma ^ \bullet _{\delta,\omega (1+ \varepsilon ) }}$
and   
is 
an unbounded subset of 
$ \gamma ^ \bullet _{\delta,\omega (1+ \varepsilon ) } $ 
and it is 
of cardinality $< \theta$,  
\sn
\item[(c)]  $\cup\{u_{\delta,\xi,\varp}:\varp < \theta\}$ is equal to 
$u_{\delta,\xi},$  
\sn
\item[(d)]  $u_{\delta,\xi,\varp}$ is computable from
  $\pr^ \mathfrak{B} 
  (\gamma^\bullet_{\delta,\varp},\xi)$
  recalling that $ \pr ^ \mathfrak{B} $ is a pairing function, 
  using as parameters $ \bar{ f }, \bar{ g } $ which were 
  fixed in $ (*)_2$.  
\end{enumerate}
\end{enumerate}
[Why? should be clear.] 

\mn
Lastly, 
\mn
\begin{enumerate}
\item[$(*)_7$]  let: 
    \begin{enumerate} 
        \item[(a)] $ \mathbf{U} = \{\bar{ u } _{\delta, \xi }
                : \delta \in S, \xi \in C_ \delta \} $ 
        \item[(b)] $\bfU_w = 
\{\bar u^{[w]}_{\delta,\xi}:\delta \in S$ and
 $\xi \in e_\delta\}$ for $w \in [\theta]^\theta$.
    \end{enumerate} 
\end{enumerate}
\mn


We shall prove that (why the $w$? for the use in the proof of part
(4) of the claim):
\mn
\begin{enumerate}
\item[$(*)_8$]  if $w \in [\theta]^\theta$ then
$\bfU_w$ witnesses $\WSGC(\lambda,\theta) \le \lambda$.
\end{enumerate}
\mn

Fix $w$ now and we shall deal with all the demands:
\mn
\begin{enumerate}
\item[$(*)_{8.1}$]  $\bfU_w$ has cardinality $\le \lambda$; in fact is
  equal to $\lambda$.
\end{enumerate}
\mn
[Why?  As $|\bfU_w| \le |\{(\delta,\xi):\delta \in S,\xi \in
e_\delta \subseteq \theta^{+} \}| \le \lambda + \theta^{+} = \lambda$.  The
other inequality is also easy as $\cup\{u_{\delta,\xi}:\delta \in
S,\xi \in e_\delta\} = \lambda$
and each $ u_{\delta ,  \xi } $ has cardinality $ \theta < \lambda $.]  
\mn
\begin{enumerate}
\item[$(*)_{8.2}$]  $\bfU_w \subseteq \bfU_{\lambda,\theta}$ is reasonable.
\end{enumerate}
\mn
[Why?  By the choices above.]
\mn
\begin{enumerate}
\item[$(*)_{8.3}$]  $\bfU_w$ semi-guess clubs.
\end{enumerate}
\mn
[Why?  Let $M$ and $A$ be as in Definition \ref{b2}(3B)(a)$'$; 
\wilog \, $ M $ expand  $ \mathfrak{B} $ and 
let 
$ M^+ $ be the expansion of $ M $ by the relation 
$ < ^{M^+}$, the order of the ordinals $ < \lambda $ 
and $ P^{M^ +}= A$,    
and let 
$E := \{\delta < \lambda:M^+ 
\rest \delta
\prec M^+ 
\}$, clearly $E$ is a club of $\lambda$.  By
the choice of $\bar C$ there is $\delta \in S$ such that $C_\delta
\subseteq E$ (hence $\delta \in  E$). 
Note that if $ \alpha \in C_ \delta $  then
$ A \cap \alpha $ is unbounded in $ \alpha $.

Now recall  
    that $ M \rest 
    \delta \prec M$,  
$\langle u_{\delta,\xi}:\xi \in e_\delta\rangle$ is
$\subseteq$-increasing continuous with union $\delta$, each
$u_{\delta,\xi}$ is of cardinality $\le \theta$ and 
$e_\delta$ is a club of $\theta^+$ 
hence $e = \{\xi \in e_\delta:M ^+
\rest u_{\delta,\xi}
\prec M^+
\}$ is a club of $\theta^+$.  
So if $ \xi \in e   $ 
    then $ A \cap u_{\delta, \xi } $ 
is unbounded in $ u_{\delta, \xi } $. 
Now
choose $\xi \in e$, so 
$\bar u
= \bar u_{\delta,\xi}$ is as required.]
\mn
\begin{enumerate}
\item[$(*)_{8.4}$]  $\bfU$ is weakly bounded.
\end{enumerate}
\mn
[Why?  Just think, recalling $(*)_1$ and Definition \ref{b2}(8), that
is, note that $\langle C_\delta \cap \alpha:\delta \in S^+\rangle$ has
cardinality $\le \theta^+$ for each $\alpha < \lambda$ 
because $ \beta 
    \in C_{\delta_1}\cap  C_{\delta _2} \Rightarrow C_{\delta _1} \cap \beta 
    = C_{\delta _2 }\cap \beta $  and $\delta > \alpha \Rightarrow \sup(C_{\delta} \cap \alpha) \in C_{\delta},$ 
anyhow  
 below we shall get more.] 
\mn
\begin{enumerate}
\item[$(*)_9$]  $\bfU$ is bounded 
hence $\CSGC(\lambda,\theta)$ holds,
  in fact:
\sn
\begin{enumerate}
\item[(a)]  if $u_1 = u_{\delta_1,\xi_1,\varp_1},u_2 =
  u_{\delta_2,\xi_2,\varp_2}$ and
  $\pr(\gamma^\bullet_{\delta_1,\varp_1},\xi_1) =
  \pr(\gamma^\bullet_{\delta_2,\varp_2},\xi_2)$ then:
\sn
\begin{enumerate}
\item[$(\alpha)$]  $\langle \gamma^\bullet_{\delta_1,\varp}:\varp \le
  \varp_1\rangle  = \langle \gamma^\bullet_{\delta_2,\varp}:\varp \le 
  \varp_2\rangle$
\sn
\item[$(\beta)$]   $u_1=u_2$ 
\end{enumerate}
\sn
\item[(b)]  $\pr(\gamma^\bullet_{\delta,\varp},\xi_1) <  
\gamma^\bullet_{\delta,\varp +1}$.
\end{enumerate}
\end{enumerate}
\mn
[Why?  Clause (a) holds by $(*)_6(d)$ and clause (b) by $(*)_1(h)$.]
\medskip

We have finished proving $ \lambda =
\CSGC( \lambda,\theta ) $,
and even $ \CSGC(\lambda, \theta, \mathfrak{B} )$, that is clause 
(A) of part (1).

\noindent
\underline{Clause (B)}: 
    Fix a stationary $ S \subseteq S^{\theta ^+}_ \theta $ which 
    belongs to $ \check{I}_ \theta [\theta ^+ 
        ]$, see Def \ref{z16}. By \ref{z19} we can  
choose $  \langle \zeta_{\xi,\varp}  : \varepsilon < \theta \rangle $ 
for $\xi \in  S $ 
such that for any such $\xi \in
S,   
\langle \zeta_{\xi,\varp}:\varp < \theta\rangle $
is increasing continuous with limit $\xi$
and $ (\zeta _{\xi _1, \varepsilon + 1} =  \zeta _{\xi _2, \varepsilon + 1})  
\wedge ( \upsilon \le \varepsilon) \Rightarrow \zeta _{\xi _1, \upsilon}
= \zeta _{\xi _2, \upsilon} $ (and more). Without loss of generality this is $\overline{C} = \langle C_{\alpha}: \alpha < \lambda$ such that if $\xi \in S, \ \varp < \theta$ then for some $\alpha < \zeta_{\xi, \varp}$   we have $\{ \zeta_{\xi, i}: i \leq \varp \} = C_{\alpha}$ and for $\overline{u} \in \bfU$ such that $\zeta = \sup(\bigcup_{\varp}u_{\varp})$ we require that $\alpha \in u_{\varp + 1}$ (and $\overline{C}$ definable in $\mathfrak{B}$).

Now in the proof of  clause (A) of 
part (1)   
    we choose $ \bar{ f }, \bar{ g } $ as in $ (*)_3$ but in addition
$ \bar{ g } = \langle  g_ \alpha : \alpha \in [ \theta, \theta ^+\rangle )$  
satisfies that: 
    if $ \alpha = \zeta _{\xi , \varepsilon },  
    \varepsilon $ 
a limit ordinal then $ g_ \alpha $ is computable from 
$ \langle g_ \beta : \beta \in \{\zeta _{\xi , \iota }: 
\iota < \varepsilon  \}  \rangle$ e.g.  as follows: for $\iota < \theta$ let $W_{\alpha}^{\iota} = \{  \gamma < \alpha: \ \text{for some} \ \beta \in \{ \zeta_{\xi, e}: e < \varepsilon \} \ \text{we have} \ g_{\beta}(\gamma) < \iota  \}$ and then define $g_{\alpha}^{\iota}$ by induction on $\iota < \theta$ such that: 

\noindent

\begin{enumerate}
    \item[$(*)$]
    \begin{enumerate}
        \item[(a)] $g_{\alpha}^{\iota}$ is a function from $W_{\alpha}^{\iota}$ onto  some ordinal $< \theta.$ 
        
        \item[(b)] $g_{\alpha}^{\iota}$ is increasing with $\iota.$
        
        \item[(c)]$g_{\alpha}^{\iota} \rest (W_{\alpha}^{\iota} \setminus \bigcup\{ W_{\alpha}^{j}: j < \iota \})$ is order preserving.
    \end{enumerate}
\end{enumerate}


Also 
we can restrict ourselves to $\xi \in S  $ 
such
that $u_{\delta,\xi}$ is closed under $\pr$.  Then we can restrict
ourselves to $(w,\delta,\xi)$ such that $\varp_1 < \varp_2 \in w
\Rightarrow \pr(\gamma^\bullet_{\delta,\varp_1},\zeta_{\xi,\varp_1})
\in u_{\delta,\xi,\varp_2}$. 


\medskip

\noindent
\underline{Clause (C)}:  Easy but we elaborate.  

We are assuming $ \lambda = \lambda ^ \theta, \theta = \cf( \theta )  $; 
so  $ \mathbf{U} =  
    \mathbf{U} _{\lambda, \theta }$ is  trivially 
a subset  of $ \mathbf{U} _{\lambda, \theta }$
of cardinality $\lambda$ and let $F$ be a one-to-one function from $\{ \overline{u} \rest \epsilon: \overline{u} \in \mathbf{U} \ \text{and} \ \epsilon < \theta  \},$ clearly exist. Let $ M $ be a model with universe
$ \lambda $   
and we have to find $ \bar{ u } $ as promised. 
Toward this we choose $ u _ \varepsilon $  by induction 
on $ \varepsilon $ as follow:

\begin{enumerate} 
    \item[(a)] $ u_ \varepsilon $ is a subset of $ \lambda $ 
        of cardinality $ < \theta $,
    \item[(b)] $ u_\varepsilon = \cl (u _\varepsilon , M  )  $   
    and has no last member,
    \item[(c)] if $ \varepsilon = \zeta + 1 $ then 
        some $ \alpha \in A \setminus  \sup( u_\zeta )$
        belongs to $ u_ \varepsilon $, 
    \item[(d)] if $ \varepsilon $ is a limit ordinal then 
        $ u_ \varepsilon = \cup \{u_\zeta : \zeta  < \varepsilon  \} $
            
    \item[(e)] If $\epsilon = \zeta +1$ then $F(\overline{u} \rest \epsilon) \in u_{\epsilon}.$        
\end{enumerate} 

There is no problem to carry the induction and 
$ \langle u_ \varepsilon : \varepsilon < \theta  \rangle $ 
is as required.

\underline{Clause (D)}  
Recall that $ {\mathscr S } \subseteq \{ w \subseteq \lambda : \otp(w)= \theta \} $ 
and more by our assumption. For each $ w \in {\mathscr S } $  let 
$ \bar{ u } _w $ be $ \langle u_{w, \varepsilon }: \varepsilon < \theta \rangle $
were $ u_{w, \varepsilon }= \{ \alpha \in W: \otp( w \cap \alpha ) 
    < \omega (1 + \varepsilon ) \}\cup \varepsilon  $. Now let $ \mathbf{U} = 
    \{ \bar{ u } _w: w \in {\mathscr S } \} $, it suffice to prove that 
    $ \mathbf{U} $ witness $ \MGC( \lambda, \theta ) = \lambda $.
    
    Clearly most demands hold:  $  \mathbf{U} $ 
    is a subset of 
    $ \mathbf{U} _{\lambda, \theta }$ of cardinality $ \lambda $, and for each 
    $ \bar{ u } \in \mathbf{U} $ the sequence 
    $ \langle sup (u_ \varepsilon ): \varepsilon < \theta \rangle $ 
    is increasing. The main point is, to be given $  M, A  $ as in 
    clause (a) of Def \ref{b2}(3) and to prove that sub-clauses $ (\alpha ),(\beta ) $ 
    there hold. 
    
    Let $ M^+ $ be an expansion of $ M $ by the order 
    $ <^{M^+} $ of the ordinals $ < \lambda $,
 $ R^{M^+}= A , \pair $  
 and let $ E = \{ \delta < \lambda : M^+ \rest \delta \prec M^+\} $, 
 clearly it is a club of $ \lambda $.    Now there is no harm 
 in replacing $ A $ by a smaller sub-set so let 
 $ A'  =  \{ \alpha \in  A : \alpha = \min (A \setminus 
 \beta ) $ for some $ \beta \in E \} $. 
 Clearly $ A' \in [\lambda ]^\lambda  $ so by the choice of $ {\mathscr S } $ 
 there is $ w \in {\mathscr S } $ such that $ w \cap A $ has cardinality $ \theta $. 
 
 Now $ \bar{ u } _ w \in \mathbf{U} $ is as required. 
 
 \underline{Clause E}:  
 
 By \cite{Sh:420} there is a stationary $ {\mathscr A} 
    \subseteq [ \lambda ] ^ \theta $ 
 of cardinality $ \lambda$, see details in the proof of
 part (3). Now for each 
 $ w \in {\mathscr S }   $ let  $  \mathscr A_ w  = \{
    v \in {\mathscr A} : w \subseteq v \}$, 
 so it is non-empty. Now for each  $w \in {\mathscr S } $  let 
 $ \langle \alpha _{w, \varepsilon }:
 \varepsilon < \theta \rangle $ list the members of 
 $ w $ in increasing  order.
 Also for each such pair $ (w, v)$ let 
 $\bar{ u } _{w, v }=  \langle  u_{w,v,\varepsilon }: 
    \varepsilon < \theta \rangle  $ 
 be such that: 
 \begin{enumerate} 
 \item[(a)] $ u_{w,v,\varepsilon }$ is a subset of $v $ 
    of cardinality $ < \theta $, 
 \item[(b)]   $ u_{w,v,\varepsilon }$ is increasing continuous 
        with $ \varepsilon $, 
 \item[(c)]  $ u_{w,v,\varepsilon }$ includes $ \{ \alpha _{w, \zeta }: 
    \zeta < \omega (1 + \varepsilon )\} $ 
    if $ \theta > {\aleph_0} $ and is $ \{\alpha _{w, \zeta }:
    \zeta < 1 +  \varepsilon \} $, 
 \item[(d)]   $ u_{w,v,\varepsilon }$ is included in $ \cup \{ \alpha _{w, \zeta }: 
    \zeta < \omega (1 + \varepsilon )  
        \} $ 
    
 \item[(e)]  $ \cup \{  u_{w,v,\varepsilon }: \sigma < \theta  \} = 
 \      v \cap \cup \{ \alpha _{w, \varepsilon }: \varepsilon < \theta \}  $ 
 \end{enumerate} 
 
 Lastly we define $ \mathbf{U} $ as the set 
 $ \{ \bar{ u } _{w, v }: w \in {\mathscr S }, v \in {\mathscr A} _ w   \} $; 
 so it suffice to prove that $\mathbf{U}  $ witnesses 
 $ \MGC_D(\lambda, \theta ) = \lambda $; this is as  
 in previous cases.
 




\noindent
2) As in \cite[Ch.III]{Sh:g}, and anyhow not used .

\medskip 

\noindent
3) By \cite{Sh:430} there is $\cS$ such that:  
\mn
\begin{enumerate}
\item[$(*)_1$]  
\begin{enumerate}
\item[(a)]  $\cS \subseteq [\lambda]^\theta$ has cardinality $\lambda$
\sn
\item[(b)]  $\cS$ is stationary, i.e. for every model $M_*$ with
  universe $\lambda$ and vocabulary $\le \theta$ there is 
$w \in \cS$ such that $M_* \rest w \prec M_*$.
\end{enumerate}
\end{enumerate}
\mn
Now as we can increase $\cS$, \wilog \,:
\mn
\begin{enumerate}
\item[$(*)_2$]  $\cS \cap [\alpha]^\theta$ is a stationary subset of
  $[\alpha]^\theta$ for every $\alpha \le \lambda$. 
\end{enumerate}
\mn

We continue as in the proof of part (1), maybe details will be given
    in \cite{Sh:F2071} and anyhow this will not be used here.  
    


\noindent
4) Let $ \mathbf{U} \subseteq \mathbf{U} _{\lambda, \theta }$    
    medium       guess clubs,
      
      
      Now clause (a) follows by \ref{b2}(3)(b).

For clause (b),  for $\overline{u} \in \mathbf{U}_{1},$ let $\langle \alpha_{u, \varepsilon}: \varepsilon < \theta \rangle$ list $\bigcup_{\varepsilon} u_{\varepsilon}$ and for each $\overline{u}$ and increasing $g \in {}^{\theta}\theta$ we define $w_{\overline{u}, \varepsilon}^{g}$ by induction on $\varepsilon < \theta$ as follows. For $\varepsilon = 0$ let $w_{\overline{u}, \varepsilon}^{g} = \cl_{\mathfrak{B}}(\emptyset),$ for $\epsilon$ limit let $w_{\overline{u}, \varepsilon}$ = $\bigcup\{  w_{\overline{u}, \zeta}: \zeta < \epsilon\},$ so let $\epsilon = \zeta +1.$ Let $\iota < \theta$ be minimal $\geq g(\zeta), \epsilon,$ such that $\sup(w_{\overline{u}, \zeta}^{g}) < \sup \{ \alpha_{\overline{u}, \iota (1)}: \iota(1) <i  \}$ and let $w_{\overline{u}, \varepsilon}^{g} = \cl_{\mathfrak{B}}(\{ \alpha_{\overline{u}, \iota(1)}: i(1) < \iota \}).$ Lastly let $\mathscr{F} \subseteq \{ g \in {^{\theta}\theta}: $g$ \ \text{is increasing} \} $ be $<_{J_{\theta}^{\mathrm{bd}}}-$unbounded of cardinality $\mathfrak{b}_{\theta}$ and let $\mathbf{U} = \{  \langle \overline{w}_{\overline{u}, \varepsilon}^{g}: \epsilon < \theta \rangle: g \in \mathscr{F} \ \text{and} \ \overline{u} \in \mathbf{U}_{1} \}. $ Now check.

    





\noindent 
5), 6)
Combine things above.  
\end{PROOF}

\begin{discussion}
\label{b20}

Assume $\lambda > \theta \ge \sigma =
\cf(\sigma),(2^\sigma > \lambda$ in the interesting case).
    Let
$\mathbf{U}  
_{\lambda,\theta,\sigma} = \{\bar u:\bar u = \langle u_\varp:\varp
< \sigma\rangle$ is $\subseteq$-increasing and $u_\varp \in
[\lambda]^\theta\}$ and repeat the definition.  Of doubtful help,
otherwise $(\theta^{++},\theta^+,\theta)$ would have helped.
\end{discussion}

\begin{theorem}
\label{b23}
1) Assume $\lambda = \cf(\lambda) > 
    \theta = \cf(\theta), \
 D= [\theta ] ^ \theta, \ $ 
  $ \AGC_D (\lambda,\theta) = \lambda$
  and $\gb_\theta > \lambda$.  
\Then \, $\lambda \notin \Univ(T_{\ceq})$; moreover,
$\univ_{T_{\ceq}}(\lambda) \ge \gb_\theta$. 

\noindent 
1A) In part (1) we can replace 
 $ \AGC_D(\lambda, \theta )$ by $ \MGC  
(\lambda, \theta )$.   

\noindent
2) If $\lambda = \cf(\lambda) > \theta = \cf(\theta) $ 
    and\footnote{
Recall that this means that $ D = \{ \theta \} $    
    }
$\FGC(\lambda,\theta) = \lambda$ and $\chi = \gd_\theta > \lambda$
or just $ \cf({}^{ \theta } \theta  
    , \le _D))\ge \chi 
    > \lambda  $,   
\then \,  
$\univ_{T_{\ceq}}(\lambda) \ge \chi $. 

\noindent
3) If $D$ is a uniform filter on $\theta,({}^\theta \theta,<_D)$ is
$(< \chi)$-directed and $\chi > \lambda,\FGC_D(\lambda,\theta) = 
\lambda$ \then \, $\univ_{T_{\ceq}}(\lambda) \ge \chi$.


\end{theorem}

\begin{remark}
\label{b26}
1) The Claim \ref{b35} below shows that we cannot weaken the  
assumption on $T$ too much.

\noindent
2) Note that the  above works also for $ \theta = {\aleph_0} $.  
   
\noindent
3) See more in \cite{Sh:F2071}.


\end{remark}

\begin{PROOF}{\ref{b23}}
1)  
So let ($T = T_{\ceq}$ and):
\mn
\begin{enumerate}

\item[$(*)_0$] $\mathfrak{B}$  is as in Definition \ref{b2}, such that: 

\begin{enumerate}
    \item[(a)] $\mathfrak{B}$ has universe $\lambda,$
    
    \item[(b)] $\tau_{\mathfrak{B}},$ the vocabulary of $\mathfrak{B},$ is countable,
    
    \item[(c)] $\mathfrak{B}$ has a pairing function $\mathrm{pr}: \lambda \times \lambda \to \lambda$ i.e., a one-to-one 2-$place$ function from $\lambda$ into $\lambda$ and $\mathrm{pr}_{1}$ and $\mathrm{pr}_{2}$ its inverses.  
\end{enumerate}

\item[$(*)_1$]  assume $\alpha_* < \gb_\theta$ and $M^*_\alpha \in
\EC_T(\lambda!)$ for $\alpha < \alpha_*$; it suffices to find $N \in
\EC_{T_{\ceq}}(\lambda!)$ not embeddable into $M^*_\alpha$ for every $\alpha <
\alpha_*$
\sn
\item[$(*)_2$]  let $\bfU \subseteq \bfU_{\lambda,\theta}$ witness
$ \AGC_D(\lambda,\theta, \mathfrak{B}) = \lambda$ or just $\MGC_{D}(\lambda, \theta, \mathfrak{B}) = \lambda$. 


\mn
[Why? Is there such $\mathbf{U}?$ by the assumption of the theorem and apply \ref{b5}(5).]
   
\sn
\item[$(*)_3$]  for $\bar u \in \bfU,\alpha < \alpha_*$ and $d \in
  P^{M^*_\alpha}$ 
  we define   
  the set 
  $E_{\bar u,d,\alpha}$; clearly 
is a club of $\theta$, as follows:  
\newline   

\begin{enumerate}
    \item[$\bullet$] $ E_{\bar u,d,\alpha} = \{\varp < \theta:\varp$ is a limit ordinal such
that $u_\varp$ is closed inside $\cup\{u_\zeta:\zeta < \theta\}$
 under the functions of $M^*_\alpha$ and the function 
$F^{M^*_\alpha}(-,d)$, or just if $a \in u_\varp,b \in
\bigcup\limits_{\zeta < \theta 
        } u_\zeta$ and $M^*_\alpha \models ``F^{M^*_\alpha }   
            (a,d) = b"$
then $b \in \sup(u_\varp)\}$. 
\end{enumerate}

\end{enumerate}
\mn
So $\cE = \{E_{\bar u,d,\alpha,  
    }:\bar u \in \bfU,\alpha < \alpha_*$
and $ d 
    \in P^{M^*_\alpha}\}$ is a set of clubs of $\theta$ of
cardinality $\le |\bfU| + |\alpha_*| +
|P^{M_{\bullet}}| < \gb_\theta  $. 
Hence, 

\mn
\begin{enumerate}
    \item[$(*)_{3.1}$] there is an increasing function $g:\theta \rightarrow \theta$ such that $(\forall E \in \cE)(\forall^\infty \varp < \theta)(g(\varp) > \suc_E(\varp))$.  
\end{enumerate}

Now we can construct $N = N_g \in \EC_T(\lambda!)$ such that:
\mn
\begin{enumerate}
\item[$(*)_4$]  
\begin{enumerate}
\item[(a)]  $P^N = \{3 \beta:\beta < \lambda\}$ hence $Q^N = \{3 \beta
  +1,3 \beta +2:\beta < \lambda\}$
\sn
\item[(b)]  if $\alpha = 3 \beta  + 1
< \lambda$ (hence $\alpha \in Q^N$)
  then $\alpha = \min(\alpha/E^N)$
\sn
\item[(c)]  if $\bar u \in \bfU$ then for some $\alpha(\bar u) =
  \alpha_{\bar u,g} \in P^N$ we have: if $\beta \in (u_{\varp +1}
  \backslash u_\varp) \cap Q^N$ and $(\beta/E^N) \cap u_\varp =
  \emptyset$ but  
    $\upsilon < \theta \Rightarrow (\beta/E^N) \cap \sup (\bigcup\limits_{\zeta  < \theta }  
  u_\zeta )  \nsubseteq \sup(u_\upsilon ) $  then $F^N  
    (\beta   , \alpha(\bar u) ) \in   
    \sup ( \bigcup\limits_{\zeta < \theta } 
  u_\zeta)\backslash \sup( u_{g(\varp + 1)})$.
\end{enumerate}
\end{enumerate}
\mn
Now toward contradiction assume that:
\mn
\begin{enumerate}
\item[$(*)_5$]  $f$ embeds $N_g$ into $M^*_\alpha$ and $\alpha <
  \alpha_*$.
\end{enumerate}
\mn
Let  $ N^+_g= (N_g, <^{N_{g}^{+}} )$  where $<^{N_{g}^{+}} \ = \{ (\alpha, \beta): \alpha < \beta < \lambda \}$ and  let  
    $M_{\bullet}$  
be a model with universe $\lambda$ expanding $M^*_\alpha$
and (a renaming of) $N^+  
        _g$;  
        (that is   $ \tau (M_{\bullet} ) $   
          contains also a 
        disjoint copy  $ \tau ' $ of $ \tau (N^+_g)$ such that  
        the restriction of  $ M_{\bullet} $   to $ \tau ' $ 
        is  the suitable copy of $ N^+_g$).
        
Also 
we have $f = G^{M_*}$ for some   
unary 
function symbol 
$G \in
\tau(M_{\bullet} )$  
and $M_{\bullet}$ has Skolem functions and $\tau(M_{\bullet})$ is countable and expand $\mathfrak{B.}$



\mn
\begin{enumerate}

\item[$(*)_{6}$] 

    \begin{enumerate}
        \item[(a)] Let $E = \{  \delta < \lambda: M_{\bullet} \rest \delta \prec M_{\bullet}  \},$ so a club of $\lambda,$ 
    
         \item[(b)] Let $H: \lambda \to \lambda$ be $H(\alpha) = \suc_{E}(\alpha),$
         
         \item[(c)] Let $M_{*} = (M_{\bullet}, H).$
    \end{enumerate}

\item[$(*)_{7}$] Let $A = \{ \mathrm{pr}^{\mathfrak{B}}(3 \beta + 1, f(3 \beta +1)): \beta < \lambda \}.$ 

\end{enumerate}

As $f$ is a  one-to-one function from $\lambda$ to $\lambda$ necessarily $A \in [\lambda]^{\lambda}.$ By the choice of $\mathbf{U}$ and of $D$ as $[\theta]^{\theta}$ there is $\overline{u}$ such that: 

\begin{enumerate}

\item[$(*)_8$]
\begin{enumerate}
\item[(a)]  $\bar u \in \bfU$,  

\sn
\item[(b)]  if $\varepsilon < \theta$ then $u_{\delta} = \cl(u_{\varepsilon}, \mathfrak{B}),$

\item[(c)]  $c \ell(u_\varp,M_{*}) \subseteq \sup(u_\varp)$, and 
    (essentially follows) $ M_{\bullet} \rest \sup(u_ \varepsilon ) \prec M_{\bullet}, $ 
\sn
\item[(d)] 
the set $ v =  
    \{\varepsilon < \theta :  
A \cap u_{\varepsilon + 1 } \setminus \sup(u_ \varepsilon )
\not= \emptyset \} $
has cardinality
  $\theta$,   

\item[(e)] $M_{*} \rest \bigcup\limits_{\varp < \theta 
    } u_\varp \prec M_{*},$ (not used when we assume only 
    $ \MGC(\lambda, \theta ) = \lambda $), 

\end{enumerate}
\end{enumerate}
\mn

Now let $ d = f( \alpha _{\bar{ u }, g})$, so it is a member of $ P^{M^*_\alpha }$, 
and 

\begin{enumerate} 
    \item[$(*)_9$] if $ \varepsilon \in v $ then for some $ a = a_\varepsilon \in 
    u_{\varepsilon + 1 } \setminus \sup(u_\varepsilon )$ 
    we have $ a _ \varepsilon  \in A$,  
\end{enumerate}

\begin{enumerate}
    \item[$(*)_{10}$] 
        \begin{enumerate}
            \item[(a)] for $\varepsilon \in v$  let $a_{\varp} =                 \mathrm{pr}^{\mathfrak{B}}(\beta_{\varepsilon}, f(\beta_{\varepsilon})),$ so $\beta_{\varp} \in \{ 3\gamma +1: \gamma  < \lambda \},$
            
            \item[(b)] Let $v_{1} = \{ \varp \in v: g(\varp) > \suc_{E(\overline{u}, d, \alpha)} (\varp) \},$
            
            \item[(c)] For $\varp < \theta$ let $\zeta_{\varp} = \suc_{E(\overline{u}, d, \alpha)}(\varp),$
        \end{enumerate}
        
    \item[$(*)_{11}$]  $v_{1} \in [v]^{\theta},$
\end{enumerate}    

[Why? Because $v \in [\theta]^{\theta}$ and the choice of $f$ (here the use of $\mathfrak{b}_{\theta}$ matters.)]
\mn


Now,

\sn
\begin{enumerate}
    \item[$\boxplus$] Assume $\epsilon \in v_{1},$
\end{enumerate}

\begin{enumerate}    
    \item[$\bullet_{1}$] $a_{\varp} \in u_{\varp +1} \setminus u_{\varp} ,$ 
\end{enumerate}

[Why? by the choice of $a_{\varp}.$]
    
\sn    
\begin{enumerate}     
    \item[$\bullet_{2}$] if $\zeta < \theta$ then $\sup(u_{\zeta})$ is closed under $x+1$ so a limit ordinal and $f^{-1};$ that is, $\forall \beta < \lambda [\beta < \sup(u_{\varp}) \iff f(\beta) < \sup(e_{\varp})],$ 
\end{enumerate}    

[Why? As $\cl(\sup(u_{\varp}), M_{*}) \subseteq \sup(u_{\varp })$ by $(*)_{8}(c)$.]
    
\sn    
\begin{enumerate}
    \item[$\bullet_{3}$] $\beta_{\varp} \in u_{\varp +1} \setminus \sup(u_{\varp})$ and $f(\beta_{\varp}) \in u_{\varp + 1} \setminus \sup(u_{\varp}),$
\end{enumerate}    
    
    [Why? As $\beta_{\varp} = \mathrm{pr}_{1}^{\mathfrak{B}}(a_{\varp}), f(\beta_{\varp}) = \mathrm{pr}_{2}^{\mathfrak{B}}(a_{\varp})$ and $a_{\varp} \in u_{\varp +1}$ by $\bullet_{1}$ (and $u_{\varp +1} = \cl_{\mathfrak{B}}(u_{\varp +1})$) by $(*)_{8}(b)$ we have $\beta_{\varp} \in u_{\varp + 1}, \ f(\beta_{\varp}) \in u_{\varp +1}.$ Now, $\beta_{\varp} < \sup(u_{\varp}) \iff f(\beta_{\varp}) < \sup(u_{\varp})$ by $\bullet_{2}$ and $(\beta_{\varp}, f(\beta_{\varp}) < \sup(u_{\varp})) \Rightarrow a_{\varp} = \mathrm{pr}^{\mathfrak{B}}(\beta_{\varp}, f(\beta_{\varp})) \ \subseteq \sup(u_{\varp + 1}).$   But, $a_{\varp} \notin  \sup(u_{\varp}),$ together $\beta_{\varp}, f(\beta_{\varp}) \notin \sup(u_{\varp})$ and we are done proving $\bullet_{3}$.]
    
    \sn
\begin{enumerate}    
    \item[$\bullet_{4}$] $F^{M_{\alpha}^{*}}(f(\beta_{\varp}), d) \notin [ \zeta_{\varp}, \sup \left( \bigcup_{\varp} u_{\varp} \right)), $
\end{enumerate}

[Why? by the definition of $E_{\overline{u}, d, \alpha}$ in $(*)_{3}$.]

\sn
\begin{enumerate}    
    \item[$\bullet_{5}$] $F^{M_{\alpha}^{*}}(f(\beta_{\varp}), d) = f(F^{N_{g}}(\beta_{\varp}, \alpha_{\overline{u}, g} )),$
\end{enumerate}

[Why? As $f$ embed $N_{g}$ into $M_{\alpha}^{*}$ and the choice of $d$ as $f(\alpha_{\overline{u}, g}.)$]

\sn
\begin{enumerate}    
    \item[$\bullet_{6}$] $F^{N_{g}}(\beta_{\varp}, \alpha_{\overline{u}, g}) \in \left[ \sup(u_{g(\varp)}), \sup( \bigcup_{\zeta}u_{\zeta})  \right),$
\end{enumerate}

[Why? By the choice of $N_{g}, \alpha_{\overline{u}, g}$ and $\bullet_{3}.$]

\sn
\begin{enumerate}    
    \item[$\bullet_{7}$] $f(F^{N_{g}}(\beta_{\varp}, \alpha_{\overline{u}, g})) \in \left[ \sup(u_{g(\varp)}), \sup(\bigcup_{\zeta}u_{\zeta})  \right),$
\end{enumerate}

[Why? By $\bullet_{2}$ and $\bullet_{6}.$]

\begin{enumerate}    
    \item[$\bullet_{8}$] $f(F^{N_{g}}(\beta_{\varp}, \alpha_{\overline{u}, g})) \in \left[ \sup(u_{\zeta(\varp)}), \sup(\bigcup_{\zeta} u_{\zeta}) \right)$ recalling $\zeta(\varp) = J_{\varp} = \suc_{E(\overline{u}, d, \alpha)}(\varp),$
\end{enumerate}

[Why? By $\bullet_{7},$  the definition of $E_{\overline{u}, d, \alpha}$ and the assumption on $\varp \in V_{1}.$]

\begin{enumerate}    
    \item[$\bullet_{9}$] $F^{M_{\alpha}^{*}}(f(\beta_{*}), d) \in \left[ \zeta_{\varp}, (\sup(\bigcup_{\zeta} u_{\zeta} ) \right),$ 
\end{enumerate}

[Why? By $\bullet_{8}$ and $\bullet_{5}.$]

\sn

But $\bullet_{4}$ and $\bullet_{9}$ are contradictory, so we are done proving part $(1).$

1A), 2), 3) Similarly to the proof of part (1)  
with some changes
 will be done in \cite{Sh:F2071}. 

\end{PROOF}

\begin{conjecture} 
\label{b29}
1) Assume $T$ (is countable complete first order) with the
$\TP_2$.  If $\lambda > \theta > \aleph_0$ are regular, $\gd_\theta >
\lambda$ and $\gd_\kappa > \FGC(\lambda,\theta)$ (maybe $\theta$
inaccessible), \then \, $\univ_T(\lambda) \ge \gd_\kappa$.

\noindent
2) Assume $T$ (is countable complete first order) non-simple.  If
$\lambda > \theta > \aleph_0$ are regular, $\gd_\theta > \lambda$ and
$\gd_\kappa > \CFGC(\lambda,\theta)$, 
\then \, $\univ_T(\lambda) \ge \gd_\kappa$.
\end{conjecture} 

\begin{remark}
\label{b32}
See hopefully  \cite{Sh:1200}, 
\cite{Sh:F2071}. 
\end{remark}

\begin{claim}
\label{b35}
Assume $\mu = \mu^{<\mu} \le \theta = \cf(\theta) < \lambda =
\cf(\lambda) < \chi = \chi^\lambda  $,  
    $ \lambda$ is strongly inaccessible Mahlo 
and for transparency 
$\GCH$ holds in the interval $[\mu,\chi)$.  For
some $\bbP$:
\mn
\begin{enumerate}
\item[(a)]  $\bbP$ is a $(< \mu)$-complete forcing of cardinality
  $\chi$ neither collapsing any cardinal, 
  nor
  changing cofinalities
\sn
\item[(b)]  $(2^\mu)^{\bfV[\bbP]} = \chi$
\sn
\item[(c)] in $ \mathbf{V} ^ \mathbb{P} $ we have $  \mathfrak{d} _ \partial  = \chi $ 
    for every inaccessible $ \partial   
            \in [\mu , \lambda )$,
\sn
\item[(d)]   in $ \mathbf{V} ^ \mathbb{P} $ there is $ \bar{ {\mathscr P} } $
        as in \ref{a14}(2),
\sn
\item[(e)] 
     $ T_{\ceq}$ has no universal member in $ \lambda $  
            moreover $ \univ( \lambda, T_{\ceq } ) \ge   \chi $ 
\item[(f)] the results of \cite{Sh:175a} holds, i.e. there is a universal
  random graph in $\lambda$,  and
 see \cite{Sh:1151}
\end{enumerate}
\end{claim}

\begin{remark} 
Note that here the case $\lqq  \cup _{\varepsilon < \theta }u_\varepsilon  \cap \sup( u_ \zeta ) $
has cardinality $ \theta $"  does not arrise.
\end{remark} 

\begin{PROOF}{\ref{b35}}
    
Our proof is based on the proof 
\cite{Sh:175a},  but the quoted \cite{B5} has to be 
changed, see full proof in a work by Mark Po\'or and the author in preparation.  

That is, we choose:  
\begin{enumerate} 
    \item[($*$)] 
$ \mathbb{P} = \mathbb{P} _3 $, where 
$ \langle \mathbb{P} _ k, \name{ \mathbb{Q} }_{\ell} :
    k \le 3, {\ell} < 3 \rangle $ is an iteration
    and:
\begin{enumerate}
\item[(A)]  $ \mathbb{Q} _0$  is adding $ \chi $ $ \lambda $-Cohen, so
    it satisfies:  
      \begin{enumerate} 
          \item[$\bullet _1$] $ \mathbb{Q} _0 $  is a $ (< \lambda )$-complete
     forcing notion of cardinality $ \chi $,
          \item[$\bullet _2$] $ \mathbb{Q} _0 $ neither  collapse some cardinal 
    nor change any cofinality  (in fact is $ \lambda ^+$-cc),
          \item[$\bullet _3$] in $ \mathbf{V} ^{\mathbb{Q} _0}$          there
    is a family $ {\mathscr A}_0 $ of $ \chi$-many 
    subsets of 
    $ \lambda  $ each of cardinality $ \lambda  $, the intersection
    of any two having cardinality $ < \lambda  $, 
      \end{enumerate} 
\item[(B)]  in $ \mathbf{V} ^{\mathbb{Q} _0}  
    = \mathbf{V} ^{\mathbb{P} _1}$ the forcing notion $ \mathbb{Q} _ 1$
    satisfies:  
          \begin{itemize}
\item[$\bullet _1$] $ \mathbb{Q} _1 $  is a $ (< \mu  )$-complete
    $  \lambda $-cc 
    forcing notion of cardinality $ \chi $, 
    (yes, $ \lambda$-cc not $ \lambda ^+$-cc),
\item[$\bullet _2$]  $ \mathbb{Q} _1 $  does  
        neither  collapses some cardinal 
        nor changes any cofinality,  
\item[$\bullet _3$]  in $ \mathbf{V} ^{\mathbb{P} _2}$  there
    is a family $ {\mathscr A}_1 $ of $ \chi$-many 
    subsets of 
    $ \lambda $, 
    each of cardinality $ \lambda $, 
    the intersection
    of any two having cardinality $ <  \mu $, 
                \item[$\bullet _4$] in $\mathbf{V} ^{\mathbb{P} _2}$ we have 
    $ \mathfrak{d} _ \partial  = \chi $ for every weakly 
    inaccessible 
    $ \partial \in (\mu, \lambda )$, 
        \end{itemize}

\item[(C)] in $ \mathbf{V} ^{\mathbb{P} _2}$ we have $ \mathbb{Q} _2 $ 
which is $ (< \mu )$-complete  $ \mu ^+ $-cc forcing notion, 
forcing that there is a universal graph of cardinality 
$ \lambda $.
\end{enumerate} 
\end{enumerate} 

Now, why are there  such $ \mathbb{Q} _ {\ell} $-s?
For clause (A)  use the forcing of adding $ \chi $ 
$ \lambda $-Cohens.
 
For clause (B) we use $ \mathbb{Q} _1 $  
    such that 
\begin{enumerate} 
    \item[$(*)$] $\mathbb{Q} = \mathbb{Q} '_1 \times  \mathbb{Q} '_1$:   
    starting with $ ( \mu, \lambda, \chi , {\mathscr A} _0)$   
        as above:
\begin{enumerate} 
\item[(a)]  $ \mathbb{Q} '_1$  is the product with Easton support
    of $ \langle  \mathbb{Q} _{1, \theta }: \theta \in S \rangle $ 
    where $ S= \{ \theta : \theta \in ( \mu , \lambda )$ is an inaccessible 
    non-Mahlo$ \} $ and for $ \theta \in S, \mathbb{Q} _{1.\theta }$   
    is the forcing adding $\chi $ many $ \theta $-Cohens,
\item[(b)]   
    $ \mathbb{Q} ''_2  $ 
 forces a refinement $ {\mathscr A} _1$ of $ {\mathscr A} _0$ 
 to a family as required as in \cite{B5}  but with Easton support; 
 so each condition has cardinality $ < \lambda $  and has Easton support. 
\end{enumerate} 
\end{enumerate} 

Now why clause (B)$\bullet _1, \bullet _2, \bullet _3 $ holds?
As said above  noting that \cite[6.1]{B5} use full support while  use Easton support

Lastly for clause (C) we apply \cite{Sh:175a}.

\medskip 

Having constructed $ \mathbb{P} = \mathbb{P} _3$ we have 
to check that is is as required.

Now being $ (< \mu ) $-complete, of cardinality $ \chi $, 
pedantically of density $ \chi $, is obvious by the properties 
of the $ \name{ \mathbb{Q} }_ {\ell} $-s. Similarly 
concerning \lqq no cardinal is collapsed and no cofinality changed",
so clause (a) of Claim \ref{b35} holds.
Also forcing the existence  of a universal graph of cardinality 
$ \lambda $ holds by the choice of $ \name{ \mathbb{Q} }_2$, 
so clause (f) of \ref{b35} holds.  

Next, clause (c)  
there saying $ \mathfrak{d} _\partial = \chi $ 
holds because    it obviously holds in 
    $ \mathbf{V} ^{\mathbb{P} _2}$   
by the choice of 
$ \name{ \mathbb{Q} }'_1 $ and the later forcing 
preserve it because it satisfies the $ \mu ^+ $-cc. 
    Now, 
    lastly, 
why  clause (d) of \ref{b35}  
holds? First, in $ \mathbf{V} ^{\mathbb{P} _1}$ we 
have  $ \GCH$   in the interval $ [\mu, \lambda)$ so there is such 
$  \bar{  {\mathscr P} }$, and   
$ \mathbb{P}_3/\mathbb{P}_1 $ satisfies the $ \lambda $-cc 
so the old clubs of $ \lambda $ are dense.
and 
this continue to holds in 
$\mathbf{V} ^\mathbb{P}  $. 
        Hence the non-existence 
of a universal model of $ T_{\ceq}$ in $ \lambda $, 
(holds by \ref{a14}(2)).
\end{PROOF}

\begin{question}
\label{b38}
1) Can we for theories $T$ satisfying $\NSOP_1 + \TP_2$ get similar results? 

\noindent
2) Is $T_{\ceq}$ in some sense minimal non-simple in a suitable family
of theories?
\end{question}

\begin{claim}
\label{b47}
 Assume $\lambda > \theta = \cf(\theta) > \aleph_0$ and $\lambda =
\lambda^\theta$ and 
$\bbP$ is a $\theta$-cc forcing notion.  

\noindent 
1) In $ \mathbf{V} ^ \mathbb{P}    $  
there is a reasonable strongly bounding 
$ \mathbf{U} \subseteq \mathbf{U} _{\lambda, \theta }$ 
of cardinality $ \lambda $  witnessing 
$ \lambda = \FGC ( \lambda, \theta , \mathfrak{B} )$  

\noindent
2) Assume $\bfU \subseteq \bfU_{\lambda,\theta}$ 
fully   
guess clubs
\then \, in
$\bfV^{\bbP},\bfU$ still 
fully  
guess clubs.

\noindent
3) In part (2), if $\bfU$ is reasonable/bounding/strongly
bounding/weakly bounding in $\bfV$, \then \, so it is in $\bfV^{\bbP}$.
\end{claim}

\begin{remark}
\label{b50}
1) This  showed help in consistency results.

\noindent
2) Similarly for the other versions of guessing clubs from \ref{b2},
but take care of what is $D$.

\noindent 
3)  In \ref{b47}(2) we can replace \lqq fully guess clubs" by a slightly stronger 
version of \lqq almost guess clubs" specifically, instead of one set $ A $ 
we have $ \sigma < \theta $  sets $ A $.

\noindent 
4) Also 
In \ref{b47}(2) we can  use versions with $ D$-s.
\end{remark}

\begin{PROOF}{\ref{b47}}  
Part (1) follows by parts (2),(3) because 
in $ \mathbf{V} $ there is such $ \mathbf{U} $  
by \ref{b11}(1)(C). 
The point is:
\mn
\begin{enumerate}
\item[$(*)$]  $(A) \Rightarrow (B)$ where:
\sn
\begin{enumerate}
\item[(A)]  if $\chi > \lambda$ and $\{\name{\gB},\bbP,\lambda\} \cup
\{\varp:\varp \le \theta\} \subseteq N \prec (\cH(\chi),\in)$ and
$\|N\| <  
    \lambda $ and $ \mathbb{P} $ satisfies the $ \theta ^+$-cc 
 where $\Vdash_{\bbP} 
``\name{\gB}$ a model with universe $\lambda$
  and vocabulary of cardinality $< \theta$"
\sn
\item[(B)]  $\Vdash_{\bbP} ``\name N \cap \lambda = \cl(|N|,\name{\gB})$ and $\name{\gB} \rest |N|$ is an elementary
  submodel of $\name{\gB}"$.
\end{enumerate}
\end{enumerate}
\end{PROOF}


\bigskip




\begin{definition} \label{c23} 
Assuming $ \theta = \cf( \theta ) \le 
    \lambda $
we let  $ \mathfrak{d} ^\dagger _{\theta \lambda }$ be the 
cofinality of the partial order $ ( {\mathscr F}^\dagger  _{\theta, \lambda } , 
\le ^\dagger _{\theta, \lambda })$  where:
\begin{enumerate} 
        \item[$(*)_1$] $ {\mathscr F}  ^\dagger _{\theta, \lambda }$  is the family
    of subsets of  $ {}^{ \theta } \theta $ of cardinality $ \le \lambda $
        \item[$(*)_2$] let $\le =  \le ^\dagger _{\theta, \lambda }$ is the following partial order
    on $ {\mathscr F}\dagger  _{\theta, \lambda } $:
    
    $ F_1 \le F_2$  if $ ( \forall f_1 \in F_1) (\exists f_2 \in F_2)
[f_1 \le f_2]$
\end{enumerate} 
\end{definition} 

The following was part of \ref{b23}, maybe we shall return to it.

\begin{claim} \label{b23a}
\noindent   
1) If $ \lambda > \theta > {\aleph_0} $ are regular, $  
\FGC(\lambda, \theta )=\lambda $  
    and $ \mathfrak{d} _ \theta > \lambda $
then $ \univ( \lambda  , 
T_{\ceq}) \ge \mathfrak{d} _ \theta  $.

\noindent 
2) Above we can replace $ \mathfrak{d} _\theta $  by 
$ \mathfrak{d} ^ \dagger _{\theta, \lambda }$
\end{claim} 

\begin{PROOF}   {\ref{b23a}}

\noindent 
1) Like the proof of part (1)  of \ref{b23}   
    so we mainly note the changes.

\begin{enumerate} 
    \item[$(*)_1$]  as above.
\end{enumerate} 

\begin{enumerate} 
    \item[$(*)_2$]  $ \mathbf{U} \subseteq \mathbf{U} _{\lambda, \theta }$ 
        witness $ \FGC(\lambda, \theta )= \lambda $.
\end{enumerate} 

\begin{enumerate} 
    \item[$(*)_3$] we let:
     \begin{enumerate}
        \item[(a)] $ M^+_ \alpha $   is an expansion of $ M^*_\alpha $
            by a pairing function and Skolem functions,
        \item[(b)] $ {\mathscr Z } = \{(\alpha, \bar{ u } , d): 
            \alpha < \alpha _*, \cup _{\varepsilon < \theta} 
            u_ \varepsilon$ is closed under the function
            $ F^{M^+_\alpha } $ and each $ u_ \varepsilon $ is closed 
            under the functions of 
            $ M^+_\alpha \} $,
        \item[(c)] for $ (\alpha, \bar{ u } , d ) \in {\mathscr Z } $ 
            let $E_{\alpha, \bar{ u } , d } = 
            \{\varepsilon < \theta : u_{\varepsilon  } $ is closed 
            under  $F^{M^*_\alpha} (-, d ) \} $,
        \item[(d)]     above let $ g_{\alpha, \bar{ u }, d}  \in {}^{ \theta } \theta $
be such that  $  g_{\alpha, \bar{ u }, d} (\varepsilon ) = 
    \min (E_{\alpha, \bar{ u } , d } \setminus (\varepsilon + 1)  ) $,   
$ g \in {}^{ \theta } \theta $  
    \end{enumerate} 
\end{enumerate} 

Next choose 
$ g \in {}^{ \theta } \theta $  not  
        bounded by any 
    well defined  $ g_{\alpha, \bar{ u }, d}$
Now we choose $ N =N_g $ as follows:
    
\begin{enumerate} 
\item[$(*)_4$]  we let
    \begin{enumerate} 
        \item[(a),(b)]  as above,
        \item[(c)] for every $ \bar{ u } \in \mathbf{U} $ for some
        $ \alpha ( \bar{ u } )= \alpha _{\bar{ u } , v} \in P^N$ 
        for every $ \varepsilon < \theta $ we have:
        \begin{itemize}
            \item $ F^N(-, d)$  maps $ \cup _{\varepsilon < \theta }u_ \varepsilon $ 
      into itself
            \item for $ \varepsilon < \theta$ we have:
            $ \varepsilon \in v $ iff 
            there is $ a \in u _ {\varepsilon + 1 }\setminus u_ \varepsilon $ 
    such that $ F^N 
                  (a, d ) \notin u_{\varepsilon + 1 } $. 
        \end{itemize}
 \end{enumerate} 
 \end{enumerate} 
    
The rest is as in the proof of part (1) of \ref{b23}.      
 \end{PROOF}

\bibliographystyle{amsalpha}
\bibliography{shlhetal}

\end{document}